\documentclass[11pt,leqno,twoside]{article}
\usepackage{fancyhdr}
\usepackage{amsmath,amssymb}
\usepackage{fouriernc}
\def\tg#1{\tag*{\textrm{#1}}} 

\let\nulset\emptyset        
\renewcommand{\emptyset}{\text{\large$\nulset$}}

\DeclareMathOperator{\reg}{reg}
\DeclareMathOperator{\sing}{sing}
\DeclareMathOperator{\dist}{dist}
\DeclareMathOperator{\graph}{graph}
\DeclareMathOperator{\spt}{spt}

\def\R{\mathbb{R}}\def\Z{\mathbb{Z}}\def\C{\mathbb{C}}\def\Sph{\mathbb{S}}
\def\S{$\mathsection$}
\renewcommand{\epsilon}{\varepsilon}
\def\wtilde#1{\hskip1pt\widetilde{\hskip-1pt#1}}
\def\wbar#1{\hskip1pt\overline{\hskip-1pt#1}}
\def\what#1{\hskip0.5pt\widehat{\hskip-0.5pt#1\hskip-0.5pt}\hskip0.5pt}
\def\tsum{{\textstyle\sum}}

\newcounter{sequation}[section] 
\renewcommand{\thesequation} {\arabic{section}.\arabic{sequation}} 
\def\tl#1{\refstepcounter{sequation}\label{#1}{\textbf{\thesequation}\hskip1.2pt }}
\def\dl#1{\refstepcounter{sequation}\label{#1}\leqno\textrm{\thesequation}} 
\def\dtg#1{\refstepcounter{sequation}\label{#1}\tag*{\textrm{\thesequation}}} 
\newcounter{pequation} 
\renewcommand{\thepequation}{(\arabic{pequation})}%
 \def\pdl#1{\refstepcounter{pequation}\label{#1}\leqno\textrm{\thepequation}}
\def\ptg#1{\refstepcounter{pequation}\label{#1}\tag*{\textrm{\thepequation}}} 

\newenvironment{state}[1]{\par\smallbreak\noindent {\bf#1} %
\textit\bgroup}{\egroup \par\ifdim\lastskip<\medskipamount \removelastskip\penalty55\medskip\fi} %
\newenvironment{proof}{\setcounter{pequation}{0}\noindent}{\,\,\nobreak$\square$\par} 

\usepackage[landscape]{geometry}%
\usepackage{2up}
\source{\magstep0}{5.5in}{8.267in}
\target{\magstep0}{11.0in}{8.5in}

\title{\vskip-.6in           Stable minimal hypersurfaces in $\R^{N+1+\ell}$ with \\ 
singular set an arbitrary closed $K\subset\{0\}\times\R^{\ell}$  \author{\scshape                         Leon Simon }
                                                \date{\vspace{-.2in}}}%

\begin{document}

\maketitle

\thispagestyle{empty}

\setlength{\abovedisplayskip}{8pt plus 3pt minus 5pt}
\setlength{\belowdisplayskip}{8pt plus 3pt minus 5pt}

\begin{abstract}
\noindent With respect to a $C^{\infty}$ metric which is close to the standard Euclidean metric on $\R^{N+1+\ell}$, where
$N\ge 7$ and $\ell\ge 1$ are given, we construct a class of embedded $(N+\ell)$-dimensional hypersurfaces (without
boundary) which are minimal and strictly stable, and which have singular set equal to an arbitrary preassigned closed subset
$K\subset\{0\}\times\R^{\ell}$.  Thus the question is settled, with a strong affirmative, as to whether there can be ``gaps'' or
even fractional dimensional parts in the singular set.  Such questions, for both stable and unstable minimal submanifolds,
remain open in all dimensions in the case of real analytic metrics and in particular for the standard Euclidean metric.

\noindent The construction used here involves the analysis of solutions $u$ of the Symmetric Minimal Surface Equation on
domains $\Omega\subset\R^{n}$ whose symmetric graphs (i.e.\ $\{(x,\xi)\in\Omega\times \R^{m}: |\xi|=u(x)\}$) lie on one
side of a cylindrical minimal cone, including in particular a Liouville type theorem for complete solutions (i.e.\ the case
$\Omega=\R^{n}$).
\end{abstract} 

\section*{Introduction}\label{intro}

With respect to a $C^{\infty}$ metric which is close to the standard Euclidean metric on $\R^{N+1+\ell}$, where
$N\ge 7$ and $\ell\ge 1$ are given, we here construct a class of embedded $(N+\ell)$-dimensional hypersurfaces
(without boundary) which are minimal and strictly stable, and which have singular set equal to an arbitrary
preassigned closed subset $K\subset\{0\}\times\R^{\ell}$.  A precise statement of the theorem is given
in~\S\ref{notation} below, and includes examples in the lowest dimension possible for embedded stable minimal
hypersurfaces with non-isolated singularities---which is dimension $8$ in $\R^{9}$. 

Thus the question is settled, with a strong affirmative, as to whether there can be ``gaps'' (as in~\cite{Sim95b}) or
even fractional dimensional parts in the singular set.  Such questions, for both stable and unstable minimal
submanifolds, remain open in all dimensions in the case of real analytic metrics and in particular for the standard
Euclidean metric.

The methods used in the present paper are primarily PDE methods, utilizing solutions and supersolutions of the
symmetric minimal surface equation (SME) and a contraction mapping argument in combination with a
Liouville-type theorem (from \cite{Sim21a}) for stable minimal hypersurfaces which lie on one side of a cylindrical
hypercone. The SME is ideal for these constructions, since it admits a rich class of singular solutions while at the
same time, as discussed in~\S{\ref{SME}}, having nice continuity and Lipschitz estimates, and it can also be
conveniently modified to handle the class of smooth ambient metrics introduced here. Additionally the method
enables us to obtain a rather precise description of the shape of the singular
examples---see~Theorem~\ref{alt-main-th} and Remark~\ref{rem-on-shape}\,(3) below.

The proof of the main theorem, including the selection of appropriate metrics but deferring the proof of strict stability, is
given in \S{\ref{outline}} below, contingent on having a suitable family of solutions of the SME.  In \S{\ref{soln-family}} the
existence of a such a family is established, using preliminaries established in~\S\S\ref{radial-solns}--\ref{small-D-y-sec}. The
strict stability of the examples obtained in~\S{\ref{outline}} will be proved in \S\ref{small-D-y-sec} (see
Remark~\ref{rem-on-shape}\,(2)).

Whether or not there can be examples like those established here in the case of low dimensional submanifolds
which are minimal with respect to smooth or real analytic metrics also remains largely an open question.  In this
direction, Zhenhua Liu~\cite{Liu20} has recently constructed examples of 3-dimensional minimizers (in higher
codimension) which have singular set consisting of the union of an arbitrary number of arcs.\footnote{ Added in
proof: \lineskiplimit-1000pt \baselineskip9pt Liu (arXiv::2110.13137) has since extended his work to include
fractional dimensional singular sets.}

\section{Notation and Statement of Main Theorem}\label{notation}

For $N\in\{1,2,\ldots\}$, $Z\in\R^{N}$ and $\rho>0$ we let
\[%
B_{\rho}^{N}(Z)=\bigl\{X\in\R^{N}:|X-Z|\le \rho\bigr\}, \,\,\,\breve{B}_{\rho}^{N}(Z)= %
\bigl\{X\in\R^{N}:|X-Z|< \rho\bigr\}, %
\]%
sometimes written $B_{\rho}(Z),\breve B_{\rho}(Z)$ when no confusion is likely to arise, and
\[%
B_{\rho}^{N}=B_{\rho}^{N}(0),\quad \breve B_{\rho}^{N}=\breve B^{N}_{\rho}(0).
\]%
$\mu_{j}$ (sometimes written $\mu$) denotes $j$-dimensional Hausdorff measure on $\smash{\R^{N}}$.

Let $M$ be a smooth embedded hypersurface in an open subset $U\subset\smash{\R^{N+1}}$, meaning that
$M\subset U$ is non-empty and for each $X\in M$ there is $\rho>0$ with $\smash{\breve
B_{\rho}^{N+1}}(X)\cap M=\psi(V)$ for some smooth, proper, rank $N$, injective map $\psi$ from an open set
$V\subset\smash{\R^{N}}$ into $\smash{\R^{N+1}}$.

For such $M$ we let $\reg M$ be the relatively open subset of $U\cap\wbar{M}$ (the closure of $M$ in $U$)
consisting of all points $X\in U\cap\wbar{M}$ such that, for some $\sigma>0$, $\breve
B_{\sigma}^{N+1}(X)\cap\wbar{M}$ is a smooth embedded hypersurface, and we let
\[%
\sing M=U\cap\wbar{M}\setminus \reg M.
\]%
We shall always assume
\[%
\reg M=M\text{ and }\sing M=U\cap\wbar{M}\setminus M,
\]%
since otherwise we could work with $\reg M$ instead of $M$.

Henceforth $n\ge 3$, $m\ge 2$, $n+m\ge 8$, $\ell\ge 1$, and points in $\R^{n}\times\R^{m}\times\R^{\ell}$ will
be denoted $(x,\xi,y)$.

\medskip

The main theorem is then as follows---a more explicit version of this theorem, with good information about the
shape of the singular examples, is given later in Theorem~\ref{alt-main-th} and Remark~\ref{rem-on-shape}.

\begin{state}{\bf{}\tl{state-main-th} Theorem.}%
Let $K$ be an arbitrary closed subset of $\,\R^{\ell}$. Then for each sufficiently small $\tau\in (0,1)$ there is a
$C^{\infty}(\R^{n+\ell})$ function $f=f(x,y)$ with $\sup|f-1|<\tau$ and $\sup_{}|D^{j}f|<C\tau\,\,\forall j\ge 1$,
$C=C(\ell,m,n,j)$, and a smooth oriented embedded hypersurface $M\subset \R^{n+m+\ell}$ which, with respect
to the metric
{\abovedisplayskip3pt\belowdisplayskip3pt%
\[%
g_{|(x,\xi,y)}=\tsum_{i=1}^{n}dx_{i}^{2}+f(x,y)\tsum_{j=1}^{m}d\xi_{j}^{2}+\tsum_{k=1}^{\ell}dy_{k}^{2}, %
\quad (x,\xi,y)\in \R^{n}\times\R^{m}\times\R^{\ell},
\]%
is minimal (i.e.\ stationary as a multiplicity~$1$ varifold in $\R^{n+m+\ell}$) and strictly stable, and and which has
\[%
\sing M =\{0\}\times\{0\}\times K.
\]}%
\end{state}%

{\bf Note:} By saying that $M$ is strictly stable we mean that there is a constant $\kappa=\kappa(M)>0$ such
that
\[%
\kappa\int_{M}(|\nabla_{M}\zeta|^{2}+\tilde r^{-2}\zeta^{2})\,d\mu(x,\xi,y) %
\le \int_{M}\bigl( |\nabla_{M}\zeta|^{2} - |A_{M}|^{2}\zeta^{2}\bigr)\,d\mu %
\dl{stab-ineq}
\]%
for all $\zeta\in C_{c}^{\infty}(\R^{n+m+\ell})$, where $\tilde r=|x|+|\xi|$, $|A_{M}|,|\nabla_{M}\zeta|$ denote
respectively the length of the second fundamental form and length of the gradient of $\zeta$ on the submanifold $M$
(relative to the Euclidean metric for $\R^{n+m+\ell}$), and $\mu$ is $(n+m-1+\ell)$-dimensional Hausdorff measure with
respect to the Euclidean metric.  The left side of~\ref{stab-ineq} is the second variation
$\frac{d^{2}}{dt^{2}}\bigl|_{t=0}\mu(M_{t})$ of the area $\mu(M_{t})$, $M_{t}=\{(x,\xi,y)+t\zeta(x,\xi,y)\nu(x,\xi,y):(x,\xi,y)\in
M\}$ (with $\nu$ a smooth unit normal for $M$); this agrees with the second variation of $M$ with respect to the metric
$g$ at least up to terms $E$ which satisfy $|E|\le C\tau\int_{M}(\tilde r^{-2}\zeta^{2}+|\nabla^{M}\zeta|^{2})\,d\mu$,
$C=C(\ell,m,n)$.  So indeed the inequality~\ref{stab-ineq} is a strict stability condition on $M$ with respect to the metric
$g$ for suitably small $\tau$.

\section{Symmetric Minimal Surface Equation (SME)} \label{SME}

For functions $u\in C^{2}(\Omega) $, where $\Omega$ is an open subset of $\R^{N}$, the mean curvature operator
$\mathcal{H}_{0}$ is defined as usual by
\[%
\mathcal{H}_{0}(u)= \sum_{i=1}^{N} D_{i}\bigl(D_{i}u/\sqrt{1+|Du|^{2}}\bigr),
\]%
and we let
\[%
\mathcal{H}(u) = \sqrt{1+|Du|^{2}}\mathcal{H}_{0}(u) \,
=\sum_{i,j=1}^{n}\bigl(\delta_{i,j}-\nu_{i}(u)\nu_{j}(u)\bigr)D_{i}D_{j}u,
\dl{mco}
\]%
where
\[%
\nu_{i}(u) = \tfrac{D_{i}u}{\sqrt{1+|Du|^{2}}}.
\dl{def-nu-i}
\]%
For positive $u\in C^{2}(\Omega)$ the Symmetric Minimal Surface Equation (SME) , is
\[%
\mathcal{M}(u)=0, %
\dl{sme}
\]%
where 
\[%
\mathcal{M}(u)=\mathcal{H}(u)-\frac{m-1}{u}. 
\dl{Mu}%
\]%
Subsequently we shall apply the discussion of this section to the case when $N=n+\ell$, so $u=u(x,y)$ with
$x\in\R^{n}$ and $y\in\R^{\ell}$.

Geometrically, the equation~\ref{sme} expresses the fact that the graph $G(u)$ of
$u$ is a hypersurface in $\R^{N+1}$ with mean curvature $(m-1)e_{N+1}\cdot \nu/u$, where
\[%
\nu=(-Du,1)/\sqrt{1+|Du|^{2}} 
\]%
is the upward pointing unit normal of $G(u)$.

More important for our present application is that the SME on a domain $\Omega\subset\R^{N}$ actually
expresses the fact that the \emph{symmetric graph} $SG(u)\subset\Omega\times\R^{m}$, defined by
\[%
SG(u)=\bigl\{(x,\xi)\in \Omega\times\R^{m}:|\xi|=u(x)\bigr\},  %
\]%
is a \emph{minimal} (i.e.\ zero mean curvature) hypersurface in $\Omega\times\R^{m}$. This is checked as
follows: Let $\tau_{1},\ldots,\tau_{N}$ be the standard orthonormal basis $e_{1},\ldots,e_{N}$ for $\R^{N}$ and
$\tau_{N+1},\ldots,\tau_{N+m-1}$ a locally defined orthonormal basis of the tangent space of $\Sph^{m-1}$, and
let $U:\Omega\times\Sph^{m-1}\to\R^{N+m}$ be defined by
\[%
U(x,\omega)=(x,u(x)\omega).
\]%
Then $U$ is $C^{\infty}$ and injective, and $U(\Omega\times\Sph^{m-1})=SG(u)$, so by the area formula
\[%
\mu_{N+m-1}(SG(u)) = \int_{\Omega}\int_{\Sph^{m-1}}\sqrt{\det P}\,d\omega dx,
\]%
where $P=(p_{ij})=(D_{\tau_{i}}U\cdot D_{\tau_{j}}U)$, so $p_{ij}=D_{x_{i}}(x,u(x)\omega)\cdot
D_{x_{j}}(x,u(x)\omega)=\delta_{ij}+D_{i}uD_{j}u$ for $i,j=1,\ldots,N$ and
$p_{ij}=D_{\tau_{i}}(x,u(x)\omega)\cdot D_{\tau_{j}}(x,u(x)\omega)=u^{2}(x)\delta_{ij}$ for
$i,j=N+1,\ldots,N+m-1$ and $p_{ij}=p_{ji}=0$ for $i=1,\ldots,N$ and $j=N+1,\ldots,N+m-1$. Hence, for $u> 0$
with $u\in C^{1}(\Omega)$,
\[%
\mu_{N+m-1}(SG(u))= \mu_{m-1}(\Sph^{m-1})\int_{\Omega}\sqrt{1+|Du|^{2}}\,u^{m-1}dx. %
\dl{vol-SG}
\]%
But on the other hand one can directly compute that the SME is the Euler-Lagrange equation for the functional
on the right and so the SME expresses the fact that $SG(u)$ is a stationary point for the area functional
$\mu_{N+m-1}(SG(u))$, and hence solutions of the SME have minimal symmetric graphs as claimed.

Being a solution of the SME is a geometrically scale invariant property: Thus if $G=\graph u$ is the graph of a solution $u$
of the SME then any homothety of $G$ is also the graph of a solution, or, equivalently, with $t\Omega=\{tx:x\in\Omega\}$,
\begin{align*}{%
&\text{If $u(x)$ satisfies the SME on $\Omega\subset\R^{N}$ and $t>0$ then}\dtg{scale-inv} \\ %
\noalign{\vskip-3pt}
&\hskip1.8in t u(t^{-1}x) \text{ also satisfies the SME on $t\Omega$.} %
}\end{align*}
If $u\ge 0$ is continuous on $\Omega$ we say that $u$ is a \emph{singular solution} of the SME on $\Omega$ if
$u^{-1}\{0\}\neq\emptyset $ and $u$ is locally the uniform limit of smooth positive solutions of the SME on
$\Omega$.

An example of a singular solution of the SME is
\[%
u(x)=\alpha_{0}|x|, \text{ where } \alpha_{0}=\sqrt{\frac{m-1}{n-1}}.
\dl{sing-sol-1}
\]%
Observe that in this case the symmetric graph $SG(u)$ is the minimal cone
\[%
\C=\{(x,\xi)\in \R^{n}\times\R^{m}:|x|^{2}/(n-1)=|\xi|^{2}/(m-1)\}.  \dl{def-C}
\]%
For a discussion of the main properties of singular and regular solutions of the SME we refer to \cite{FouS20}. 
The main results in~\cite{FouS20} include a gradient estimate for both singular and regular solutions, but here we
shall only need the more standard gradient estimate from~\cite[Theorem~1]{Sim76}, which includes
(see~\cite[Example 4.1]{Sim76}) the result that if $u$ is a $C^{2}(B^{N}_{\rho})$ solution of the prescribed mean
curvature equation
\[%
\tsum_{i=1}^{N}D_{i}(D_{i}u/\sqrt{1+|Du|^{2}})=H,
\]%
where $|H|\le b/\sqrt{1+|Du|^{2}}$ and $|u|\le M$ on $B^{N}_{\rho}$, then for each $\theta\in [\frac{1}{2},1)$ $|Du|$ is
bounded in $B_{\theta\rho}$ in terms of $N$, $\theta$, $\rho b$ and $M/\rho$. In particular this applies to the SME on the
ball $B^{n+\ell}_{\rho}$ provided there are constants $M> L>0$ with $L\le u\le M$, in which case we have the above
hypotheses with $b=(m-1)/L$, so for $\theta\in [\frac{1}{2},1)$, 
\[%
\sup_{B_{\theta\rho}}|Du|\le C,\quad C=C(\ell,m,n,\theta,M/\rho,L/\rho).  %
\dl{grad-est}
\]%
(Note in particular that no smoothness or continuity properties of $H$ are needed for this bound.)

If $u_{1}, u_{2}$ are $C^{2}$  on a domain $\Omega\subset\R^{N}$ and $\mathcal{H}$ is
as in~\ref{mco}, then
\[%
\mathcal{H}(u_{1})-\mathcal{H}(u_{2}) = \mathcal{L}_{u_{1},u_{2}}(u_{1}-u_{2}),
\dl{H-diff}
\]%
where 
\[%
\mathcal{L}_{u_{1},u_{2}}(v)=\sum_{i,j}\bigl(\delta_{ij}-\nu_{i}(u_{1})\nu_{j}(u_{1})\bigr)D_{i}D_{j}v+ %
\sum_{k}b_{k}D_{k}v,  
\dl{op-for-diff}
\]%
with
\[%
\begin{aligned}%
&b_{k}=\sum_{i,j} \bigl(\nu_{j}( u_{1})b_{ik}(u_{1},u_{2})+\nu_{i}( u_{2})b_{jk}(u_{1},u_{2})\bigr)D_{i}D_{j}u_{2}, \\
& b_{ik} =-\int_{0}^{1}D_{p_{k}}\bigl(p_{i}/\sqrt{1+|p|^{2}}\bigr)\bigl|_{p=Du_{1}+tD(u_{2}-u_{1})}\,dt.
\end{aligned}%
\]%
Hence, if $u_{1}, u_{2}$ are positive,
\[%
\mathcal{M}(u_{1})-\mathcal{M}(u_{2}) %
=\mathcal{L}_{u_{1},u_{2}}(u_{1}-u_{2})+ \frac{m-1}{u_{1}u_{2}}(u_{1}-u_{2}). %
\dl{eqn-for-diff}
\]%
In particular if $u_{1},u_{2}$ are solutions of the SME then
\[%
\mathcal{L}_{u_{1},u_{2}}(u_{1}-u_{2})=- \frac{m-1}{u_{1}u_{2}}(u_{1}-u_{2}). %
\dl{lin-eqn-for-diff}
\]%
Also, if $u_{1}, u_{2}$ are in $C^{2}(\wbar{\Omega})$ and $\mathcal{M}(u_{1})\le\mathcal{M}(u_{2})$,
then
\[%
\mathcal{L}_{u_{1},u_{2}}(u_{1}-u_{2})\le -\frac{m-1}{u_{1}u_{2}}(u_{1}-u_{2}), %
\dl{lin-diff}
\]%
where $\mathcal{L}v$ is as in~\ref{op-for-diff}, so by the classical maximum principle $u_{1}-u_{2}$ cannot have a
zero minimum in $\Omega$ unless $u_{1}=u_{2}$ in $\Omega$, because \ref{lin-diff} says
$\mathcal{L}(u_{1}-u_{2})\le 0$ in $\Omega$ in case $u_{1}\ge u_{2}$.

We also need to discuss second variation of the symmetric area functional $\mathcal{F}(u)
=\int_{\Omega}\sqrt{1+|Du|^{2}}\,u^{m-1}dx$. By definition of $\mathcal{M}(u)$, the first variation
$\frac{d}{dt}\mathcal{F}(u+t\zeta)|_{t=0}$, assuming we are looking at positive functions $u\in C^{2}(\Omega)$
with $\Omega\subset\R^{N}$, is given by
\[%
\frac{d}{dt}\mathcal{F}(u+t\zeta)\bigl|_{t=0}=-\int_{\Omega}V^{-1}\mathcal{M}(u) \,\zeta\,u^{m-1}dx,\quad
\zeta\in C^{1}_{c}(\Omega),
\]%
where $\mathcal{M}(u)$ as in~\ref{Mu}, and $V=\sqrt{1+|Du|^{2}}$.  If $\mathcal{M}(u)=0$ then we can compute the
second variation
\[%
\frac{d^{2}}{dt^{2}}{\mathcal{F}}(u+t\zeta)\bigl|_{t=0}=-\int_{\Omega}\zeta\mathcal{L}_{u}(\zeta)\,u^{m-1}Vdx, %
\dl{2nd-var}
\]%
where $\mathcal{L}_{u}(\psi) = V^{-1}\frac{d}{dt}\bigl|_{t=0}\mathcal{M}(u+t\psi)$, and, after some calculation
and rearrangement of terms, this can be written  
{\abovedisplayskip6pt\belowdisplayskip3pt%
\begin{align*}{%
&\mathcal{L}_{u}(\psi) = (u^{m-1}V)^{-1}\tsum_{i,j}D_{i}\bigl(Vu^{m-1}g^{ij}D_{j}(V^{-1}\psi)\bigr) %
\dtg{jac-0} \\  %
\noalign{\vskip-4pt}
&\hskip1.9in +\Bigl(|A_{G(u)}|^{2}+\bigl(\frac{m-1}{V^{2}u^{2}}\bigr)\Bigr)(V^{-1}\psi), %
}\end{align*}%
with $g^{ij}=\delta_{ij}-\nu_{i}\nu_{j}$, $\nu_{i}=V^{-1}D_{i}u$, and
\[%
|A_{G(u)}|^{2}=V^{-2}\tsum_{i,j,p,q}g^{ij}g^{pq}u_{ip}u_{jq},\,\,\, u_{ij}=D_{i}D_{j}u,
\]%
is the squared length of the second fundamental form of
\[%
G(u)=\graph u=\bigl\{(x,z)\in\R^{N}\times\R:z=u(x)\bigr\}.
\]%
Notice that the equation~\ref{jac-0} can be thought of as a linear operator applied to $V^{-1}\psi$ (rather than to $\psi$), and
in that case the coefficient of the degree zero term is $|A_{G}(u)|^{2}+ \smash{\frac{m-1}{V^{2}u^{2}}}$, which one can check
is just the squared length
\[%
|A_{SG(u)}|^{2} = |A_{G}(u)|^{2}+ \frac{m-1}{V^{2}u^{2}}
\dl{2nd-ff-s}
\]%
of the second fundamental form of the symmetric graph $SG(u)$.  Also the
remaining terms (i.e.\ the first and second order terms) are in fact just  the Laplace-Beltrami
operator $\Delta_{SG(u)}(V^{-1}\psi)$ of the symmetric graph $SG(u)$, written in terms of the local coordinates
$x\in \Omega$ (and valid for functions $\psi$ which are also written in terms of the local variables
$x\in\Omega$).

So~\ref{jac-0} can alternatively be written 
\[%
\mathcal{L}_{u}(\psi) =\Delta_{SG}(V^{-1}\psi) +|A_{SG}(u)|^{2}(V^{-1}\psi).  %
\dl{jac-1}
\]%

We shall need the following consequence of the maximum principles of Ilmanen~\cite{{Ilm96}} and
Solomon/White~\cite{SolW89}.

\begin{state}{\bf{}\tl{il-sw-app} Lemma.} %
Let $\theta\in [\tfrac{1}{2},1)$ and $\delta\in(0,1]$ be given. There exist positive $\eta=\eta(\ell,m,n,\theta,\delta),
\,\tilde\eta=\tilde\eta(\ell,m,n,\theta, \delta)\in(0,\delta]$ such that if $u$ is a solution of the \emph{SME} with
$u(x,y)-\alpha_{0}|x|>0$ on $\breve B^{n+\ell}_{1}$, then
\begin{align*}{%
u(0,0)\le \eta \Rightarrow  u(x,y)-\alpha_{0}|x|\le \delta,\quad (x,y)\in B_{\theta}, %
\tg{\rm (i)} \\ %
u(0,0)\ge \delta \Rightarrow  u(x,y)-\alpha_{0}|x|\ge \tilde\eta,\quad (x,y)\in B_{\theta}. %
\tg{\rm (ii)}  %
}\end{align*}%
\end{state}%

{\bf{}\tl{grad-bd} Remark:} If $u$ is a solution of the SME with $u-\alpha_{0}|x|>0$ on $\breve B_{1}$, if
$u(0,0)\in(0,\eta_{1}]$ where $\eta_{1}$ is the constant $\eta(\ell,m,n,\theta,\delta)$ above with $\delta=1$, and if
$\tilde\eta_{0}$ is the constant $\tilde\eta(\ell,m,n,\theta,\delta)$ in case $\delta=u(0,0)$, then (i) and (ii) above give
$\tilde\eta_{0}\le u(x,y)-\alpha_{0}|x|\le 1$ for all $(x,y)\in B_{\theta}$. Hence by~\ref{grad-est}
\[%
u(0,0)\in (0,\eta_{1}] \Rightarrow{\sup}_{B_{\theta^{2}}}|Du|\le C, \quad (x,y)\in B_{\theta^{2}}, %
\,\,C=C(\ell,m,n,\theta,u(0,0)).
\]%

\smallskip

\begin{proof}{\bf{}Proof of~\ref{il-sw-app}: }  
If (i) fails for some given $\theta,\delta$ then there is a sequence $u_{k}$ of solutions of the SME on $\breve B_{1}$ with
$u_{k}-\alpha_{0}|x|>0$ and $u_{k}(0,0)\to 0$, yet such that
\[%
u_{k}(x_{k},y_{k})-\alpha_{0}|x_{k}|> \delta %
\pdl{il-sw-1}
\]%
for some $(x_{k},y_{k})\in B_{\theta}$.  

Let $M_{k}=SG(u_{k})$.  According to \cite[Lemma~2.3]{FouS20}, $\mu(M_{k}\cap \{(x,\xi,y):(x,y)\in
B^{n+\ell}_{\rho},\,|\xi|<R \})$ bounded independent of $k$ for each $R>0,\rho\in(0,1)$, so by the Allard compactness
theorem there is a subsequence of $k$ (still denoted $k$) such that $M_{k}$ converges in the varifold sense in
$\{(x,\xi,y):(x,y)\in \breve B^{n+\ell}_{1},\,\xi\in\R^{m}\}$ to a stationary integer multiplicity varifold $V$ with support of
$V$ equal to a closed subset $M$ of $\{(x,\xi,y):(x,y)\in \breve B^{n+\ell}_{1},\,\xi\in\R^{m},\,|\xi|\ge \alpha_{0}|x|\}$ and with
density of $V\ge 1$ at each point of $M$. Also $M_{k}$ converges to $\spt V$ locally in the Hausdorff distance sense in
$\{(x,\xi,y):(x,y)\in \breve B^{n+\ell}_{1},\,\xi\in\R^{m}\}$. In particular $(0,0)\in M$.

With $\C$ as in~\ref{def-C}, by virtue of the maximum principle of Solomon and White \cite{SolW89}, we then have either
$M\cap \C=\emptyset $ or $\C_{1}\subset M$, where $\C_{1}=\C\cap\{(x,\xi,y):(x,y)\in \breve B_{1}\}$.  $M\cap
\C=\emptyset$ implies $(0,0)\in M\cap \wbar{\C}\subset\{0\}\times\R^{\ell}$, which contradicts the maximum principle of
Ilmanen \cite{Ilm96}, so the alternative $\C_{1}\subset M$ must hold.  But with $(x_{k},y_{k})$ as in \ref{il-sw-1}, continuity
of $u_{k}$ implies
  \[%
[0,\delta]\subset \{u_{k}(t(x_{k},y_{k}))-\alpha_{0}|tx_{k}|:t\in [0,1]\},\quad k=1,2,\ldots, %
\pdl{il-sw-2}
\]%
so $M$ contains points $(x_{t},\xi_{t},y_{t})$ with with $(x_{t},y_{t})\in B_{\theta}$ and $|\xi_{t}|=\alpha_{0}|x_{t}|+t\delta$
for each $t\in [0,1]$. With $v(\C_{1})$ the varifold obtained by taking $\C_{1}$ with multiplicity~1, we then have $\wtilde
V=V-v(\C_{1})$ is stationary and also contains the points $(x_{t},\xi_{t},y_{t})$, $t\in(0,1]$, in its support. So the same
argument applies to $\wtilde V$ to show $\C_{1}\subset\spt \wtilde V$ and $\what V=\wtilde V-v(\C_{1})$ is stationary 
and again has the points $(x_{t},\xi_{t},y_{t})$, $t\in(0,1]$ in its support.  After finitely many such steps this gives a
contradiction to the fact that $V$ has bounded density on $\C\cap \{(x,\xi,y):(x,y)\in B_{\theta}\}$. So~\ref{il-sw-1} is
impossible and~(i) is proved.

Notice that we similarly get a contradiction if $u_{k}(0,0)\ge \delta$ and $u_{k}(x_{k},y_{k})-\alpha_{0}|x_{k}|\to
0$ instead of $u_{k}(0,0)\to 0$ and  $u_{k}(x_{k},y_{k})-\alpha_{0}|x_{k}|\ge \delta$, so~(ii) is also proved.
\end{proof}

\section[Proof of the Main Theorem] {Proof of the Main Theorem }\label{outline}

With $K$ an arbitrary closed non-empty subset of $\R^{\ell}$, let $\tau\in (0,\frac{1}{4}]$ (to be chosen later, depending
only on $\ell,m,n$) and let $h\in C^{\infty}(\R^{\ell})$ satisfy
\[%
\left\{\begin{aligned}
& h>0 \text{ on }U=\R^{\ell}\setminus K,\,\,  h=0  \text{ on }   K  \\  %
& h(y)+|D_{\!y}h(y)|+|D^{2}_{y} h(y)|+|D^{3}_{y} h(y)| < \tau, \quad y\in \R^{\ell},  \\ %
& \dist^{-j}(\partial U,y)|D^{k}h(y)|\le C\tau, \quad y\in\R^{\ell},\, j,k=0,1,2,\ldots,\, C=C(j,k). %
\end{aligned}\right.%
\dl{def-h}
\]%
It is of course standard that such functions $h$ exist.

For the proof of the main theorem (Theorem~\ref{state-main-th}) we shall need the following, which guarantees, for
each $\tau$ sufficiently small, the existence of a positive smooth solution $u_{\tau}(r,y)$ ($r=|x|$) of the SME on
\[%
\Omega=\bigl\{(x,y)\in\R^{n}\times\R^{\ell}:y\in U,\,\,|x|< h^{2}(y)\bigr\}, \quad U=\R^{n}\setminus K, %
\dl{def-om}
\]%
with $u_{\tau}(r,y)-\alpha_{0}r>0$ ($\alpha_{0}=\smash{\bigl(\frac{m-1}{n-1}\bigr)^{1/2}}$ as in~\ref{sing-sol-1}) on
$\wbar{\Omega}\setminus (\{0\}\times K)$ and $u_{\tau}-\alpha_{0}r$ vanishing to infinite order on approach to
$(0,y)\in\{0\}\times \partial U$, and with $|D_{\!y}u_{\tau}|$ small.

\begin{state}{\bf{}\tl{ex-thm} Theorem.} %
Let $\delta>0$.  There is $\tau_{0}=\tau_{0}(\delta,\ell,m,n)\in(0,\frac{1}{2}]$ such that if $\tau\in (0,\tau_{0}]$, if $h$ as
in~\emph{\ref{def-h}}, and if $\Omega$ as in~\emph{\ref{def-om}}, then there is a $C^{\infty}(\wbar{\Omega})$ solution
$u_{ \tau}$ of the \emph{SME} with
\[%
\left\{\begin{aligned}%
& |D_{\!y}u_{\tau}(r,y)|<\delta,\,\,|Du_{\tau}(r,y)|\le 2\alpha_{0},  \quad y\in U,\,r<h^{2}(y), \\ %
&\alpha_{0}r<u_{\tau}(r,y)< \alpha_{0}r+C  h^{j}(y), \,\,\,    y\in U,\,r<h^{2}(y),\,\,j\ge 1,\,\,C=C(j,\ell,m,n),
 \end{aligned}\right.%
\leqno{(\ddag)}
\]%
and such that $M=SG(u_{\tau})$ satisfies the strict stability inequality
\[%
\kappa\int_{M}(|\nabla_{M}\zeta|^{2}+\tilde r^{-2}\zeta^{2})\,d\mu(x,\xi,y) %
\le \int_{M}\bigl( |\nabla_{M}\zeta|^{2} - |A_{M}|^{2}\zeta^{2}\bigr)\,d\mu %
\]%
for all $\zeta=\zeta(x,y)\in C^{1}_{c}(\Omega)$, where $\tilde r=|\xi|+|x|$ and $\kappa=\kappa(\ell,m,n)>0$.
\end{state}%

\smallskip

{\bf{}\tl{ex-rem} Remark:} $u_{\tau}-\alpha_{0}r$ satisfies the linear elliptic equation~\ref{lin-eqn-for-diff} (with
$u_{\tau},\alpha_{0}r$ in place of $u_{1},u_{2}$) on $\{(x,y)\in\R^{n}\times\R^{\ell}:0< |x|<h^{2}(y)\}$, so by
using~$(\ddag)$ together with interior quasilinear estimates for solutions of the SME with bounded gradient in combination
with standard estimates for~\ref{lin-eqn-for-diff}, we have
\[%
\sup_{\frac{1}{4}h^{2}(y)< r< \frac{3}{4} h^{2}(y)}  %
               (h(y)))^{-j}\bigl|D^{k}\bigl(u_{\tau}(r,y)-\alpha_{0}|x|\bigr)\bigr|  %
               \le  C\tau, \forall j,k=1,2,\ldots, 
\]%
where $C=C(\ell,m,n,j,k)$.

\smallskip

The proof of Theorem~\ref{ex-thm} will be given in~\S{\ref{soln-family}}.

\medskip

Assuming the result of Theorem~\ref{ex-thm} for the moment, we now give the proof of the main theorem,
Theorem~\ref{state-main-th}.  We are going to construct the hypersurface $M$ of Theorem~\ref{state-main-th}
as the symmetric graph of a Lipschitz function $u$ on $\R^{n+\ell}$. We in fact can immediately use
Theorem~\ref{ex-thm} to define a suitable $u$:

Let $\zeta:\R\to[0,1]$ be a $C^{\infty}$ function with $\zeta(t)=1$ for $t\le\tfrac{1}{2}$, $\zeta(t)=0$ for $t\ge\frac{3}{4}$,
$|D_{t}^{k}\zeta(t)|\le 20$ for $k=1,2$, and  $|D_{t}^{k}\zeta(t)|\le C_{k},k\ge 3$.  With $\tau\in (0,\tau_{0}]$
and $u_{\tau}$ as in Theorem~\ref{ex-thm} and Remark~\ref{ex-rem}, define
\[%
u(r,y) = \begin{cases}%
{}\zeta(r/h^{2}(y))u_{\tau}(r,y)+ (1-\zeta(r/h^{2}(y)))\alpha_{0}r, &(r,y)\in [0,\infty)\times U \\
{}\alpha_{0}r, &(r,y)\in [0,\infty) \times K,
\end{cases}
\dl{def-u}
\]%
where  $h$ as in~\ref{def-h}.  Then
\[%
u(r,y)-\alpha_{0}r=\zeta(r/h^{2}(y))(u_{\tau}(r,y)-\alpha_{0}r), \quad (r,y)\in [0,\infty)\times U,
\]%
and, by~\ref{ex-thm} and~\ref{ex-rem}, $\mathcal{M}(u)=0$ for  $r\in[0,\tfrac{1}{2}h^{2}(y))$, 
$u(r,y)=\alpha_{0}r$ for $r\in [\tfrac{3}{4}h^{2}(y),\infty)$, and
\[%
\begin{cases}
0<u(r,y)-\alpha_{0} r \le C h^{j} & r\in [0,\infty) \,\,\forall j\ge 1,\,\,C=C(j) \\ %
|D^{k}(u(r,y)-\alpha_{0} r)|\le C h^{j}(y),  &0<\tfrac{1}{4}h^{2}(y)<r,\, \forall j,k,\,\,C=C(j,k). %
\end{cases}%
\dl{props-v}
\]%

We shall actually prove a more precise version of Theorem~\ref{state-main-th}, as follows.

\begin{state}{\bf{}\tl{alt-main-th} Theorem (Main Theorem).}%
There is $\tau_{0}=\tau_{0}(\ell,m,n)\in (0,\frac{1}{2}]$ such that if $\tau\in (0,\tau_{0}]$, if $h$ is as in~\emph{\ref{def-h}},
and if  $u$ is as in~\emph{\ref{def-u}}, then, with $M=SG(u)$,  there is a $C^{\infty}([0,\infty)\times\R^{\ell})$
function $f=f(r,y)$ with 
\[%
\begin{cases}%
\quad \qquad f(r,y)=1, &(r,y)\in ([0,\infty)\times K)\cup([0,\tfrac{1}{2}h^{2}(y)]\times U) \\  %
\sup |f-1|<\tau, \,\,\,\sup|D^{j}f|<C\tau,   &j=1,2,\ldots, \,\,C=C(\ell,m,n,j),  %
\end{cases}%
\]%
and, with respect to the metric
\[%
g_{|(x,\xi,y)}=\tsum_{i=1}^{n}dx_{i}^{2}+f(x,y)\tsum_{j=1}^{m}d\xi_{j}^{2}+\tsum_{k=1}^{\ell}dy_{k}^{2}, %
\quad (x,\xi,y)\in \R^{n}\times\R^{m}\times\R^{\ell},
\]%
$M$ is minimal (i.e.\ stationary as a multiplicity~$1$ varifold in $\R^{n+m+\ell}$),  strictly stable as
in~\emph{\ref{stab-ineq}}, and $\sing M=\{0\}\times\{0\}\times K$.
\end{state}%


\begin{proof}{\bf{}Proof:} By the definition~\ref{def-u} of $u$, $M=SG(u)$ is smooth near each point of $M\setminus
(\{0\}\times\{0\}\times K)$ and has conical singularities at each point of $\{0\}\times\{0\}\times K$, so $\sing
M=\{0\}\times\{0\}\times K$ by construction, and thus to complete the proof we just have to find $g$ so that $M$ is
minimal and strictly stable with respect to $g$.

To construct such a $g$, first take any positive $f=f(r,y)\in
C^{\infty}([0,\infty)\times\R^{\ell})$ and define a smooth metric
\[%
g=\sum_{i=1}^{n}dx_{i}^{2}+f(r,y)\sum_{j=1}^{m}d\xi_{j}^{2}+\sum_{k=1}^{\ell}dy_{k}^{2} %
\pdl{metric}
\]%
on $\R^{n+m+\ell} =\R^{n}\times\R^{m}\times\R^{\ell}$, $r=|x|$.

Applying the area formula as in the discussion of \S{\ref{SME}} with $N=n+\ell$, except that now we use the metric $g$ for
$\R^{n+m+\ell}$ rather than the standard metric (so now we use $g(U_{\tau_{i}},U_{\tau_{j}})$ in that argument in place of
$U_{\tau_{i}}\cdot U_{\tau_{j}}$), we have
\[%
\mu_{g}(SG(u))=\mu_{m-1}(S^{m-1})\int_{\Omega}\sqrt{1+f\,|Du|^{2}}\,\,f^{(m-1)/2}u^{m-1}\,dxdy %
\pdl{g-vol}
\]%
for any positive $C^{2}$ function $u$ on a domain $\Omega\subset\R^{n+\ell}$, where $\mu_{g}$ denotes
$(n+m+\ell-1)$-dimensional Hausdorff measure on $\R^{n+m+\ell}$ with respect to the metric $g$.  Thus the
Euler-Lagrange equation for the functional
\[%
\int_{\Omega}\sqrt{1+f|Du|^{2}}\,\,f^{(m-1)/2}u^{m-1}\,dxdy
\]%
is equivalent to the statement that the symmetric graph $SG(u)$ is a minimal (zero mean curvature) hypersurface
relative to the metric $g$ for $\R^{n+m+\ell}$.  By direct computation, the Euler-Lagrange equation is in fact
\[%
\tfrac{1}{2}\bigl(m+(1+f\,|Du|^{2})^{-1}\bigr)Df\cdot Du = %
-f\tsum_{i,j=1}^{n+\ell}\bigl(\delta_{ij}-\tfrac{fD_{i}uD_{j}u}{1+f|Du|^{2}}\bigr)D_{i}D_{j}u+\frac{m-1}{u}, %
\]%
where we use the notation $(x,y)=(x_{1},\ldots,x_{n+\ell})$ (i.e.\ $x_{n+j}=y_{j}$).  In case $u=u(r,y)$, which we
assume below, we can take $f=f(r,y)$ with equation
\[%
\begin{aligned}%
&\tfrac{1}{2}\bigl(m+(1+f\,|Du|^{2})^{-1}\bigr)Du\cdot Df = %
-f\Bigl(\Delta u -f\tfrac{Q(u)}{1+f|Du|^{2}}\Bigr)+\tfrac{m-1}{u},
\end{aligned}%
\pdl{f-eqn}%
\]%
where $D=(D_{r},D_{y})=(D_{r},D_{y_{1}},\ldots,D_{y_{\ell}})$, $Du\cdot Df=D_{r}u\,D_{r}f+D_{y}u\cdot D_{y}f$ and
\[%
Q(u) =u^{2}_{r}u_{rr}+
\tsum_{i,j=1}^{\ell}u_{y_{i}}u_{y_{j}}u_{y_{i}y_{j}}+2u_{r}\tsum_{j=1}^{\ell}u_{y_{j}}u_{ry_{j}}.
\]%
\ref{f-eqn} is a non-degenerate quasilinear first order PDE for the function~$f$ at points where $Du=(D_{r}u,D_{y}u)\neq
0$, and if $f$ is a local solution of the equation in a ball $B_{\sigma}(0,y_{0})$, then, with
$f_{\sigma}(r,y)=f((0,y_{0})+(\sigma r,\sigma y))$ for $|(r,y)|<1$ (i.e.\ translation of $y$ and scaling of $(r,y)$),
\begin{align*}{%
&\text{$f_{\sigma}$ satisfies~\ref{f-eqn} on the unit ball   $B_{1}(0)$ provided we replace $u$ by } %
\ptg{scaling-prop}\\  %
\noalign{\vskip-3pt}
&\hskip0.6in \text{ the geometrically rescaled function $\sigma^{-1}u((0,y_{0})+(\sigma r,\sigma y))$}. 
}\end{align*}
This scaling behavior is of course to be expected, given the geometric context leading to~\ref{f-eqn}.

With $z=1-f$, \ref{f-eqn} can be written in the form
\[%
\mathcal{A}(u,z)\,  Du\cdot Dz =\mathcal{M}(u) -z\mathcal{E}(u,z), %
\pdl{mod-f-eqn}%
\]%
where $\mathcal{M}$ is as in~\ref{Mu} and
\[%
\begin{aligned}%
\mathcal{A}(u,z)&= \tfrac{1}{2}\bigl(m+\bigl(1+|Du|^{2}-z\,|Du|^{2}\bigr)^{-1}\bigr)  \\
\mathcal{E}(u,z)&= \Delta u-\frac{\bigl(1+(1-z)(1+|Du|^{2})\bigr) Q(u)} {(1+|Du|^{2}-z|Du|^{2})(1+|Du|^{2})}. %
\end{aligned}%
\]%
Also, since 
\[%
\mathcal{M}(\alpha_{0} r)=0, \mathcal{A}(\alpha_{0}r,z)=\tfrac{1}{2}\bigl(m+(1+\alpha_{0}^{2}-
\alpha_{0}^{2}z)^{-1}\bigr), \text{ and }\mathcal{E}(\alpha_{0} r,z)=(n-1)\alpha_{0}/r, 
\]%
after some rearrangement of the terms, \ref{mod-f-eqn} can be written in the form
\begin{align*}{%
&  \tfrac{1}{2}\bigl(m+\bigl(1+\alpha_{0}^{2}-\alpha_{0}^{2}z\bigr)^{-1}\bigr)\alpha_{0} z_{r}+a(r,y,z)\cdot Dz  %
  \ptg{z-equn} \\ %
&\hskip1.8in  = -(n-1)\alpha_{0}z/r+ zb(r,y,z) +c(r,y) %
}\end{align*}%
where
\[%
\begin{aligned}%
a(r,y,z )&= \mathcal{A}(u,z)Du-\mathcal{A}(\alpha_{0} r,z)D(\alpha_{0} r), \\ %
b(r,y,z)&= -\bigl(\mathcal{E}(u,z)-\mathcal{E}(\alpha_{0} r,z)\bigr), \\ %
c(r,y)&=\mathcal{M}(u)=\mathcal{M}(u)-\mathcal{M}(\alpha_{0} r),
\end{aligned}
\]%
and so by~\ref{props-v}, for $j,k=0,1,2,\ldots$,
\[%
\sup_{\{(r,z):\frac{1}{4} h^{2}(y)<r<h^{2}(y),\,|z|<\frac{1}{2}\}} %
\hskip-20pt (r^{k}|D_{r,y,z}^{k}a|+r^{k+1}|D_{r,y,z}^{k}b|+r^{k+1}|D_{r,y}^{k}c|) \le C h^{j}(y), %
\pdl{z-equn-2}
\]%
where $C=C(j,k,\ell,m,n)$.  In particular there is $\tau_{0}=\tau_{0}(\ell,m,n)$ such that
\[%
|a|<\tfrac{\alpha_{0}}{2} \text{ provided }\tau\in (0,\tau_{0}] \text{ for  }|z|<\tfrac{1}{2}, %
\pdl{z-equn-3}
\]%
so the equation remains non-degenerate as long as $|z|<\frac{1}{2}$.

We first aim to get local solutions of~\ref{z-equn} with $z$ having initial value $0$ on the hypersurface
$\Sigma=\{(\frac{1}{4}h^{2}(y),y):y\in U\}$.  Thus we impose the initial condition
\[%
f(r,y)=1\text{ for } y\in U,\,\, r=\tfrac{1}{4}h^{2}(y).
\pdl{init-cond}
\]%
In view of~\ref{scaling-prop}, it is convenient to discuss this in a rescaled setting.  In fact, for given $y_{0}\in U$,
we take the translation/rescaling $(r,y)\to \rho^{-1}(r,y-y_{0})$ with $\rho=\frac{1}{4}h^{2}(y_{0})$, and in the
rescaled setting we claim, with $\tau_{0}=\tau_{0}(\ell,m,n)>0$ sufficiently small (and independent of $y_{0}$)
and assuming also~\ref{z-equn-2}, that we can find a $C^{\infty}$ solution of the local initial value problem
\[%
\left\{\begin{aligned}%
&\tfrac{1}{2}\bigl(m+\bigl(1+\alpha_{0}^{2}-\alpha_{0}^{2}z\bigr)^{-1}\bigr)\alpha_{0} z_{r}+a(r,y,z)\cdot Dz   \\  %
\noalign{\vskip-1pt}
&\hskip0.2in = -(n-1)\alpha_{0}z/r+ zb(r,y,z) +c(r,y), \,\, \Psi(y)\le r< 4\Psi(y), \,\, |y|\le 4,  \\  %
&\,\, z(\Psi(y),y) =0, \quad |y|\le 4,  %
\end{aligned}\right.%
\pdl{z-prob}
\]%
where $\Psi(y)=\frac{1}{4}h^{2}(y_{0}+\rho y)/\rho$ with $\rho=\frac{1}{4}h^{2}(y_{0})$, so that by~\ref{def-h}
\[%
\Psi(0)=1,\,\, {\sup}_{|y|<5}|D^{k}\Psi(y)|\le C\tau,\,\,\,k=1,2,\ldots, \,\,\,C=C(\ell,m,n,k). %
\pdl{P-k}
\]%

Recall that the Lagrange procedure (``method of characteristics'') guarantees local solvability in $C^{\infty}$ of
first order equations in $\R^{N}$ of the form 
\[%
{\tsum}_{i=1}^{N}a_{i}(x,z)D_{i}z=c(x,z),\quad a_{i},c\in C^{\infty},
\]%
with zero initial data on the hypersurface $\Sigma$:
\[%
\Sigma=\{(\Psi(y),y):y\in V\},   %
\]%
where $V$ is open in $\R^{N-1}$, $\Psi\in C^{\infty}(V)$, and $a(\Psi(\eta),\eta)\cdot (1,-D\Psi(\eta))\neq 0$.

Notice that geometrically this latter condition requires $a$ to not be tangent to $\Sigma$ at each point of
$\Sigma$.

The method involves first solving the ODE system
\[%
\left\{\begin{aligned}
&\, \tfrac{\partial}{\partial t} X(t,\eta) = a(X(t,\eta),Z(t,\eta)) \\  %
&\, \tfrac{\partial}{\partial t} Z(t,\eta) = c(X(t,\eta),Z(t,\eta)),   %
\end{aligned}\right.%
\]%
subject to the initial condition
\[%
X(0,\eta)=(\Psi(\eta),\eta),\,\,\,\,Z(0,\eta)=0, \,\,\,\eta\in V.
\]%
Then one proves that for each $\eta_{0}\in V$, and suitable $\rho=\rho(\eta_{0},a_{i},c)>0$, the map
$X:(t,\eta)\in\breve [0,\rho]\times B^{N-1}_{\rho}(\eta_{0})\mapsto X(t,\eta)\in R^{N}$ is a diffeomorphism onto
some open neighborhood $W$ of $(\Psi(\eta_{0}),\eta_{0})$ in $\R^{N}$, and then $z$ is defined in $W$ by
$z=Z\circ X^{-1}$.  One can then check that $z$ satisfies the PDE in $W$ with $z=0$ on $W\cap\Sigma$.

In the present case~\ref{z-prob}, we have $N=1+\ell$ and $X=(R,Y)$, with points in $\R^{1+\ell}$ denoted $(r,y)$,
$r>0$, and $\Psi(\eta)$ as in~\ref{P-k}, and the ODE system is
\[%
\begin{cases}%
\,\frac{\partial}{\partial t} R(t,\eta) =  %
\,\tfrac{1}{2}\bigl(m+(1+\alpha_{0}^{2}-\alpha_{0}^{2}Z)^{-1}\bigr)\alpha_{0}+ a_{1}(R,Y,Z)&   \\   %
\,\frac{\partial}{\partial t} Y(t,\eta) = \tilde a(R,Y,Z)\,\quad (\tilde a=(a_{2},\ldots,a_{\ell+1}))&  \\ %
\,\frac{\partial}{\partial t} Z(t,\eta) =   %
\,\bigl(-(n-1)\alpha_{0} R^{-1} +b(R,Y,Z) \bigr)Z + c(R,Y),& %
\end{cases} %
\pdl{ODE-sys}
\]%
subject to the initial conditions
\[%
R(0,\eta)=\Psi(\eta),\,\,\,Y(0,\eta)=\eta, \,\,\, Z(0,\eta) =0, \,\,\,|\eta|<5. %
\]%
We first claim that  for $\tau\le \tau_{0}=\tau_{0}(\ell,m,n,P)$ small enough, then 
the solution $(R(t,\eta),Y(t,\eta),Z(t,\eta))$ exists for $(t,\eta)\in [0,5]\times \breve B_{5}^{\ell}$.  To prove this claim, first
note that by~\ref{z-equn-2} the equation for $R$ ensures that $D_{t}R>0$ and then the initial condition for $R$
ensures that
\[%
R(t,\eta)\ge \Psi(\eta)\,\,\,(\,>1-C\tau>\tfrac{1}{2})\text{ for }(t,\eta)\in [0,5]\times B_{5}^{\ell}, %
\pdl{R-inc}
\]%
provided $\tau_{0}=\tau_{0}(\ell,m,n)$ is small enough. Then the equation for $Z$, together
with~\ref{z-equn-2}, says $|D_{t}Z|\le 2n|Z|+\tau\le 2n(|Z|+\tau)$, and hence $e^{-2nt}(|Z|+\tau)$ is
decreasing, so
\[%
|Z(t,\eta)| \le C\tau, \,\,\, (t,\eta)\in [0,5]\times \breve B_{5}^{\ell}.  \pdl{Z-bd-1}
\]%

Then by differentiating the equation for $Z$ with respect to $\eta_{j}$, integrating with respect to $t$ and using
the initial condition $Z(0,\eta)=0$ (hence $D_{\eta}Z(0,\eta)=0$) we see that also
\[%
|D_{\eta}Z(t,\eta)| \le C\tau.
\]%
So now by using the equations for $(R,Y)$ directly
\[%
|D_{t}(R,Y)-(c_{0},0,\ldots,0)|\le C\tau,\,\,\, c_{0}=\tfrac{1}{2}\alpha_{0}\bigl(m+(1+\alpha_{0}^{2})^{-1}\bigr),
\]%
and by integrating with respect to $t$,
\[%
(R,Y)(t,\eta)=(c_{0}t+\Psi(\eta),\eta)+E(t,\eta), %
\pdl{R-Y-bd}
\]%
where $E(0,\eta)=0$ and $|E|+|D_{t}E|\le C\tau$.  

Also by first differentiating the $(R,Y)$ equations with respect $\eta_{j}$ and then integrating with respect to $t$, we prove
that $|D_{\!\eta}E|<C\tau$, so in fact $|E|+|D_{t,\eta}E|<C\tau$.  So \ref{R-Y-bd} shows that
\[%
(R,Y)(t,\eta)=(c_{0}t,\eta)+\widetilde{E}(t,\eta), %
\pdl{R-Y-bd-2}
\]%
with $|\widetilde{E}(t,\eta)|+|D_{t,\eta}\widetilde{E}(t,\eta)|\le C\tau$, so, with $\tau\in (0,\tau_{0}]$, 
$\tau_{0}=\tau_{0}(m,n,\ell)>0$ small enough, $(R,Y)$ is a $C^{1}$ diffeomorphism
\[%
\Phi:[0,5]\times B^{\ell}_{5}\to W\supset \{(r,y): y\in B^{\ell}_{4},\,\, \Psi(y)\le r \le 4\Psi(y) \}, \pdl{P-5-W}
\]%
and hence $z=Z\circ\Phi^{-1}\bigl| \{(r,y): y\in B^{\ell}_{4},\,\, \Psi(y)\le r \le 4\}$ is the required solution
of~\ref{z-prob} on $\bigl\{(r,y): y\in B^{\ell}_{4}(y_{0}),\,\, \Psi(y)\le r \le 4\Psi \bigr\}$ with $z=0$ on the
hypersurface $\bigl\{( \Psi(y),y): y\in B^{\ell}_{4}\bigr\}$. 

Also, because $\mathcal{M}(u)=0$ in $\{(r,y):y\in U\text{ and }\frac{1}{4}h^{2}(y)<r<\frac{1}{2} h^{2}(y)\}$, this solution $z$
vanishes identically in the region $\bigl\{(r,y):|y|<4,\,\,\Psi(y)\le r\le 2\Psi(y)\bigr\}$ by the ODE uniqueness theorem.

Next note that, with
\[%
\mathcal{X}_{k}=D_{\eta}^{k}(R,Y,Z)\,\,\, \text{ (and $\mathcal{X}_{0}=(R,Y,Z)$)},
\]%
we can successively differentiate in~\ref{ODE-sys} to give
\[%
D_{t}\mathcal{X}_{k} = F_{k}(t,\eta) + G_{k}(t,\eta)\mathcal{X}_{k}, \pdl{X-k-equn}
\]%
for $k\ge 1$ with $F_{k},G_{k}$ smooth functions and
\[%
|G_{k}|\le C_{0},\,\,C_{0}=C_{0}(\ell,m,n),\,\,|F_{k}|\le C, \,\,\, C=C(\ell,m,n,k),
\]%
where the second inequality is subject to the inductive assumption that for $k\ge 1$ we already have bounds
$|\mathcal{X}_{j}|\le C_{k}$ for $j=0,\ldots,k-1$.  Then by subdividing the interval
\[%
[0,5]=\cup_{j=1}^{5N}[(j-1)/N,j/N],
\]%
and by integration in~\ref{X-k-equn} with respect to $t\in [(j-1)/N,s+(j-1)/N]$, where $s\in [0,1/N]$, we obtain
\[%
\begin{aligned}%
&{\sup}_{(t,\eta)\in[(j-1)/N,j/N]\times B^{\ell}_{5}} |\mathcal{X}_{k}(t,\eta)| %
\le {\sup}_{t=(j-1)/N,\,\eta\in B_{5}^{\ell}} |\mathcal{X}_{k}(t,\eta)| \\ %
\noalign{\vskip-3pt}  %
&\hskip1.7in +C + N^{-1}C_{0}{\sup}_{(t,\eta)\in[(j-1)/N,j/N]\times B^{\ell}_{5}}|\mathcal{X}_{k}(t,\eta)|, %
\end{aligned}%
\]%
where $C=C(\ell,m,n,k)$. Hence choosing $N=N(\ell,m,n)>2C_{0}$ we have
\[%
{\sup}_{(t,\eta)\in[(j-1)/N,j/N]\times B^{\ell}_{5}}|\mathcal{X}_{k}(t,\eta)| %
\le 2{\sup}_{t=(j-1)/N,\,\eta\in B_{5}^{\ell}} |\mathcal{X}_{k}(t,\eta)|+2C.  %
\pdl{X-j-ineq}
\]%
In case $j=1$ we can use the initial data $\mathcal{X}_{0}(0,\eta)=(\Psi(\eta),\eta)$, and so~\ref{X-j-ineq} gives
\[%
\sup_{(t,\eta)\in[0,1/N]\times B^{\ell}_{5}}|\mathcal{X}_{k}(t,\eta)| \le C,\,\,\, C=C(k,\ell,m,n).  %
\pdl{X-1-ineq}
\]%
For $j\ge 2$ and with $N=N(k,\ell,m,n)>2C_{0}$, \ref{X-j-ineq} gives
\[%
\begin{aligned}%
&{\sup}_{(t,\eta)\in[(j-1)/N,j/N]\times B^{\ell}_{5} }|\mathcal{X}_{k}(t,\eta)| \\ %
\noalign{\vskip-3pt} %
&\hskip1in \le 2{\sup}_{(t,\eta)\in[(j-2)/N,(j-1)/N]\times B^{\ell}_{5}}|\mathcal{X}_{k}(t,\eta)| +2C
\end{aligned}%
\]%
and so
\[%
{\sup}_{(t,\eta)\in[0,5]\times B^{\ell}_{5}}|\mathcal{X}_{k}(t,\eta)| %
\le C,\,\,\, C=C(k,\ell,m,n). %
\pdl{X-eta-ests}
\]%
Now it follows that, for $,j,k=0,1,2,\ldots$,
{\abovedisplayskip3pt\belowdisplayskip3pt%
\[%
{\sup}_{(t,\eta)\in[0,5]\times B^{\ell}_{5}}| D^{j}_{t}D_{\eta }^{k} \mathcal{X}(t,\eta)| \le C, %
\,\,\, C=C(\ell,m,n,j,k),  %
\pdl{X-ests}
\]}%
because $j=0$ holds by~\ref{X-eta-ests}, and then the case $j=1$ of~\ref{X-ests} is true by~\ref{X-k-equn}, and
finally the case $j\ge 2$ of~\ref{X-ests} is proved by induction on $j$ by applying $D_{t}^{j-1}$ to each side
of~\ref{X-k-equn}.  So~\ref{X-ests} is proved for all $j,k$.

Thus $\Phi=(R,Y)$ in~\ref{P-5-W} is actually a $C^{\infty}$ diffeomorphism with
\[%
|D_{r,y}^{k}\Phi^{-1}(r,y)|\le C,\quad y\in B_{4}^{\ell},\,\,\,\,\Psi(y)\le r\le 4\Psi(y),\,\,C=C(\ell,m,n,k),  %
\pdl{P-k-bd}
\]%
and in particular $C$ does not depend on $y_{0}$.

In view of~\ref{X-ests}, with $\mathcal{Z}_{k}(t,\eta)=D_{\eta}^{k}Z(t,\eta)$ (and $\mathcal{Z}_{0}=Z$), we can
take $k$ derivatives with respect to the $\eta$ variables in the equation for $Z$ to give
\[%
D_{t}\mathcal{Z}_{k} = F_{k}(t,\eta) +D_{\eta}^{k}\bigl(c(R(t,\eta),Y(t,\eta))\bigr),
\]%
for $k\ge 1$, where $|F_{k}|\le C\sum_{j=0}^{k}|D_{t,\eta}^{j}(c(R(t,\eta),Y(t,\eta)))|$ subject to the inductive
assumption, $|\mathcal{Z}_{j}|\le C\sum_{i=0}^{j}\bigl|D_{t,\eta}^{i}(c(R(t,\eta),Y(t,\eta)))\bigr|$ $\forall j\in\{
0,\ldots,k-1\}$, and then arguing inductively as in the proof of~\ref{X-eta-ests},~\ref{X-ests} (except that here the
argument is slightly simpler because $\mathcal{Z}_{k}$ has initial data zero by virtue of the fact that $Z(t,\eta)=0$
for all sufficiently small $t$, because $c(r,y)=0$ for $\Psi(y)\le r\le 2\Psi(y)$) to give
\[%
\sup_{(t,\eta)\in [0,5]\times B^{\ell}_{5}}\hskip-10pt |D_{t}^{j}D^{k}_{\eta}\mathcal{Z}| %
\le C\hskip-15pt\sup_{(t,\eta)\in [0,5]\times %
  B^{\ell}_{5}}\hskip-10pt\tsum_{i=0}^{j+k}|D_{t,\eta}^{i}\bigl(c(R(t,\eta),Y(t,\eta))\bigr)| %
\le C\tau h^{i}(y_{0}) %
\pdl{Z-k-bd}
\]%
by~\ref{z-equn-2}, $i,j,k=0,1,2,\ldots$, where $C=C(\ell,m,n,i,j,k)$ and in particular $C$ does not depend on
$y_{0}$; notice that here we used
\[%
\tfrac{1}{2}h^{2}(y_{0})\le h^{2}(y)\le 2h^{2}(y_{0}) \text{ for } |y-y_{0}|\le 2h^{2}(y_{0})
\pdl{y-y-0}
\]%
(for $\tau$ small enough), because by 1-variable calculus $|h^{2}(y)-h^{2}(y_{0})|=|(y-y_{0})\cdot
\int_{0}^{1}D(h^{2})(y_{0}+ty)\,dt|\le C\tau h^{2}(y_{0})$ if $|y-y_{0}|\le 2h^{2}(y_{0})$.

Thus $z=Z\circ\Phi^{-1}$ is the required solution of~\ref{z-prob} on $\bigl\{(r,y):\Psi(y)\le r\le 4\Psi(y),\,\,|y|\le 4\bigr\}$, so
changing the scale back to the original (i.e.\ $(r,y)\to (0,y_{0})+\rho(r,y)$ with $\rho=\frac{1}{4}h^{2}(y_{0})$), and using the
uniqueness theorem for solutions of the initial value problem for first order quasilinear PDE, we finally have a smooth
solution $z$ of~\ref{mod-f-eqn} on $\{(r,y):y\in U,\,\,\frac{1}{4}h^{2}(y)\le r\le h^{2}(y)\}$ with $z$ identically zero on
$\frac{1}{4}h^{2}(y)\le r\le \frac{1}{2}h^{2}(y)$.  Also, by~\ref{P-k-bd}, \ref{Z-k-bd} and~\ref{y-y-0},  $z$ satisfies
\[%
|D_{r,y}^{k}z(r,y)| \le C\tau h^{j}(y), \,\,\,\, C=C(j,k,\ell,m,n), \,\,j,k=0,1,2,\ldots, %
\pdl{sup-bd-z}
\]%
for all $y\in U$ and $\frac{1}{4}h^{2}(y)\le r\le h^{2}(y)$.

For $r\ge h^{2}(y)$ (where $u(x,y)=\alpha_{0}r$) the equation~\ref{z-equn} is just the ODE
\[%
\bigl(m+\bigl(1+\alpha_{0}^{2}-\alpha_{0}^{2}z\bigr)^{-1}\bigr)z_{r} = -2 z(n-1)/r. %
\pdl{f-ODE}
\]%
The appropriate integration shows this is equivalent to 
\[%
\tfrac{d}{dr}\bigl((1+\alpha_{0}^{2}-\alpha_{0}^{2} z(r,y))^{-\beta_{1}}\, z(r,y)\,r^{\beta_{2}}\bigr)=0 
\]%
for $r\ge h^{2}(y)$, where $\beta_{1}=\frac{1}{m(1+\alpha_{0}^{2})+1}<\frac{1}{2}$ and $\beta_{2}=
\frac{2(n-1)}{m+(1+\alpha_{0}^{2})^{-1}}$.  So in particular
\begin{align*}{%
&\bigl(1+\alpha_{0}^{2}-\alpha_{0}^{2}  z(r,y)\bigr)^{-\beta_{1}}z(r,y)=   %
\ptg{ex-form-f}\\ %
\noalign{\vskip-3pt}
&\hskip0.3in \bigl(1+\alpha_{0}^{2}-\alpha_{0}^{2} z(h^{2}(y),y)\bigr)^{-\beta_{1}} z(h^{2}(y),y)\, %
                                  \bigl(h^{2}(y)/r\bigr)^{\beta_{2}},\,\,\,r\ge h^{2}(y). %
}\end{align*}%
Thus $f=1-z$ is defined and smooth on the entire region $\{(r,y):y\in U,\,r\ge \frac{1}{4}h^{2}(y)\}$ with the
required properties, including~\ref{sup-bd-z} and the fact that $f$ is identically $1$ in the region
$\frac{1}{4}h^{2}(y)\le r\le \frac{1}{2}h^{2}(y),\,y\in U$. Finally $f=1$ evidently satisfies the equation in the region
$0\le r\le\frac{1}{4}h^{2}(y)$ (because $\mathcal{M}(u )=0$ in this region), so the proof is complete by extending $f$ to be
$1$ in this region. \end{proof} 

\vskip5pt
This completes the proof of the main theorem~\ref{state-main-th}, except for the proof of the existence result of
Theorem~\ref{ex-thm}, which will be given in~\S{\ref{soln-family}}, and the proof of the strict stability of $M$,
which will be established in Remark~\ref{rem-on-shape}\,(2) in~\S{\ref{small-D-y-sec}}.

\section{Radially Symmetric Solutions of the SME}\label{radial-solns}

To facilitate the construction of a suitable family of solutions of the SME of the type specified in
Theorem~\ref{ex-thm} of the previous section, we first need to consider the special solutions $u(x,y) =\varphi(r)$
($r=|x|$)---i.e.\ solutions of the SME which are expressible as a function of the variable $r=|x|$, or in other words
solutions $\varphi(r)$ which satisfy the Euler-Lagrange equation of the area functional
\[%
\mathcal{F}(u)=\int_{0}^{1}\sqrt{1+(u'(r))^{2}}\,u^{m-1}r^{n-1}\,dr. %
\dl{m-n-fnl}
\]%
In this case the SME is the ODE
\[%
\frac{1}{r^{n-1}}\frac{d}{dr}\bigl(r^{n-1}\frac{\varphi'}{\sqrt{1+(\varphi')^{2}}} \bigr)
=\frac{m-1}{\varphi\sqrt{1+(\varphi')^{2}}}.  \dl{varphi-equn-0}
\]%
which, for solutions with $\varphi'(0)=0,\, \varphi(0)=1$ (which we are mainly interested in here), can be
integrated to give
\[%
\frac{\varphi'}{\sqrt{1+(\varphi')^{2}}} =\frac{1}{r^{n-1}}
\int_{0}^{r}\frac{m-1}{\varphi(s)\sqrt{1+(\varphi'(s))^{2}}}\,s^{n-1}ds.
\]%
So such solutions have $\varphi'(r)\ge 0$ (and hence $\varphi(r)\ge 1$) on the entire interval $[0,r_{0}),\, (r_{0}\in
(0,\infty])$ of their existence. We can also directly prove local existence of such solutions by the contraction
mapping principle, by checking that operator
\[%
T(\varphi)= 1+\int_{0}^{r}\sqrt{1+(\varphi'(\tau))^{2}}\,\tau^{1-n}\!\int_{0}^{\tau}
\frac{m-1}{\varphi(s)\sqrt{1+(\varphi'(s))^{2}}}\,s^{n-1}dsd\tau
\]%
is a contraction mapping of the closed set $\bigl\{\varphi\in
C^{1}([0,t]):\varphi(0)=1,\varphi'(0)=0,0\le\varphi'(r)\le 1\,\forall r\in[0,t]\bigr\}$ into itself for small enough
$t>0$.

To proceed further, we write the equation~\ref{varphi-equn-0} in the form
\[%
\varphi'' +\tfrac{(n-1)}{r}(1+(\varphi')^{2})\varphi' = \tfrac{(m-1)}{\varphi}(1+(\varphi')^{2}). %
\dl{varphi-equn}
\]%
As observed in~\ref{sing-sol-1}, one solution of~\ref{varphi-equn}, although singular at $r=0$, is
\[%
\varphi_{0}=\alpha_{0} r, \quad \alpha_{0}=\sqrt{\frac{m-1}{n-1}}, %
\dl{sing-sol}
\]%
and in this case  the symmetric graph $SG(\varphi_{0})$ is just the minimal cone
$\bigl\{(x,\xi)\in\R^{n}\times\R^{m}:(n-1)|\xi|^{2}=(m-1)|x|^{2}\bigr\}$. We use the notation
\[%
\C_{0} = SG(\varphi_{0}),\,\,\C=\C_{0}\times\R^{\ell}. %
\dl{cones-C}
\]%
Notice that the solution $\varphi_{0}$ has an isolated singularity when viewed as a function of $x\in\R^{n}$, but
as a function of $(x,y)\in\R^{n+\ell}$ the singular set is the entire subspace $\{0\}\times\R^{\ell}$.

We know from the above discussion and general ODE theory that there is a unique $C^{\infty}$ solution $\varphi$
of~\ref{varphi-equn} subject to the initial conditions
\[%
\lim_{r\downarrow0}\varphi(r)=1,\quad\lim_{r\downarrow0}\varphi'(r)=0 %
\dl{initial-conds}
\]%
on a maximal interval $(0,r_{0})$, where $0<r_{0}\le \infty$ and, as we already observed above, this solution has
$\varphi'\ge 0$ and $\varphi\ge 1$ on the entire interval $(0,r_{0})$.  Differentiation gives
\[%
\varphi'''+(1+3(\varphi')^{2})\frac{n-1}{r}\varphi'' %
\ge(1+(\varphi')^{2}) \bigl((n-1)r^{-2}-(m-1)\varphi^{-2}\bigr)\varphi' >0 %
\]%
at points $r$ where $\varphi>\alpha_{0} r$, which says
\[%
(r^{n-1}e^{A(r)}\varphi'')'>0,\text{ where }A(r)=3(n-1)\int_{1}^{r}(\varphi'(t))^{2}t^{-1}\,dt
\]%
at such points.  So $r^{n-1}e^{A(r)}\varphi''$ is strictly increasing at points $r$ where $\varphi>\alpha_{0} r$ and
in particular $\varphi''>0$ on any interval $(0,\rho)$ where $\varphi>\alpha_{0}r$; also $\varphi''(0)>0$ directly from the
equation, because $(n-1)\varphi'(r)/r=(n-1)\varphi''(\theta r)$ for some $\theta\in (0,1)$ by the mean value theorem of
1-variable calculus.

Since $\frac{n-1}{r}(\alpha_{0}r)'=\frac{m-1}{\alpha_{0} r}$,  the equation for $\varphi$ can be written
\[%
(\varphi-\alpha_{0}r)'' +(1+(\varphi')^{2})\frac{(n-1)}{r}(\varphi-\alpha_{0}r)' =(m-1)(1+(\varphi')^{2})
\bigl(\frac{1}{\varphi}-\frac{1}{\alpha_{0}r}\bigr),
\]%
which is
\[%
(\varphi-\alpha_{0}r)'' +(1+(\varphi')^{2})\frac{(n-1)}{r}(\varphi-\alpha_{0}r)'  %
=-\frac{m-1}{\alpha_{0}r\,\varphi}(1+(\varphi')^{2}) \bigl(\varphi-\alpha_{0}r\bigr), %
\dl{lin-equn}
\]%
hence
\[%
(e^{A(r)/3}r^{n-1}(\varphi-\alpha_{0}r)')' <0
\]%
at points where $\varphi>\alpha_{0}r$, so $e^{A(r)/3}r^{n-1}(\varphi-\alpha_{0}r)'$ is strictly decreasing, hence $<0$
since it vanishes as $r\downarrow 0$, on any interval $(0,\rho)$ where $\varphi>\alpha_{0}r$.  In particular
\[%
\varphi' < \alpha_{0}, \text{ and hence } \varphi(r)<\alpha_{0}r +1
\]%
on any interval $(0,\rho)$ where $\varphi >\alpha_{0}r$. Thus on any such interval $(0,\rho)$ we have
\[%
\varphi''(r)>0,\, 0<\varphi'(r)<\alpha_{0}, \,\alpha_{0}r<\varphi(r)<\alpha_{0}r+1, %
\text{ and }\varphi(r)-r\varphi'(r) >0. %
\dl{varphi-asymp}
\]%
Now according to~\cite[Theorem 2.1]{HarS85} there is a smooth complete area minimizing hypersurface
$S\subset U_{+}=\bigl\{(x,\xi)\in\R^{n}\times\R^{m}:|\xi|>\alpha_{0}|x|\bigr\}$ with $\dist(0,S)=1$, and the
homotheties $\{tS\}_{|_{t>0}}$ foliate all of $U_{+}$.  Then if $\varphi(\rho)=\alpha_{0}\rho$ for some $\rho\in
(0,r_{0})$ we could choose a homothety $tS$ of $S$ which lies on one side of $SG(\varphi|[0,\rho])$ and makes
contact at some point in $U_{+}$, which contradicts the maximum principle.  So in fact $\varphi(r)>\alpha_{0}r$
for all $r\in (0,r_{0})$ and \ref{varphi-asymp} holds on the whole maximal interval $(0,r_{0})$ and in particular
$1<\varphi(r)<\alpha_{0} r+1$ and $0<\varphi'(r)<\alpha_{0}$ on $(0,r_{0})$. So $r_{0}=\infty$ by the ODE
extension theorem.

Now, since $\varphi-r\varphi'>0$, we see that $\graph \varphi$ intersects every ray $\bigl\{t ,s t):t >0\bigr\}$ with
$s>\alpha_{0}$ transversely in a single point, and so the homotheties $tSG(\varphi)$ ($=SG(\varphi_{t})$, where
$\varphi_{t}(r)=t\varphi(r/t)$) foliate all of $U_{+}$. Thus $SG(\varphi)$ is minimizing, and by construction
$\dist(0,SG(\varphi))=1$, hence the uniqueness part of~\cite[Theorem 2.1]{HarS85} is applicable, giving
$SG(\varphi)=S$.

Also, the calibration argument of Lawson \cite{Law72} shows that $\varphi_{0}=\alpha_{0} r$ \emph{strictly}
minimizes the area functional~\ref{m-n-fnl} in the sense that there is a fixed constant $C>0$ such that
\[%
\mathcal{F}(u) \ge \mathcal{F}(\varphi_{0})+ C\rho^{n+m-1}
\]%
whenever $\rho\in (0,\tfrac{1}{2}]$ and $u:[0,1]\to [0,\infty)$ is $C^{1}$ with $u(r)-\alpha_{0} r\ge 0$,
$u(r)>\rho$ for each $r\in (0,1)$, and $(u(r)-\alpha_{0} r)\bigl|_{r=1}=0$. Hence~\cite[Theorem 3.2]{HarS85} is
applicable, giving $\varphi(r)-\alpha_{0}r\sim \kappa r^{\gamma}$ as $r\to \infty$ for some $\kappa>0$, where
\[%
\gamma= -\bigl(n+m-3\bigr)/2+\sqrt{\bigl(\bigl(n+m-3\bigr)/2\bigr)^{2}-(n+m-2)}. %
\dl{def-g-+}
\]%
Thus, using~\ref{varphi-asymp},
\[%
\left\{\begin{aligned}%
&\begin{aligned}&\varphi''(r) >0,\,\,\, \varphi(r)-r\varphi'(r)>0,\,\,\alpha_{0} r<\varphi(r)<1+\alpha_{0} r, \\ %
\noalign{\vskip-4pt} %
&\hskip1.8in \text{ and } \,\,0<\varphi'(r)<\alpha_{0}\,\,\text{ for all } r>0,
\end{aligned} \\ %
&\begin{aligned}%
&\varphi(r)-\alpha_{0} r\sim \kappa r^{\gamma}, \,\,   %
                        \varphi(r)-r\varphi'(r)\sim \kappa(1-\gamma) r^{\gamma}, \\ %
\noalign{\vskip-4pt} %
& \hskip1.8in \text{ and }\varphi''(r)\sim \kappa\gamma(\gamma-1)r^{\gamma-2} \text{ as $r\to\infty$}, %
\end{aligned}%
\end{aligned}\right.%
\dl{props-varphi}
\]%
where $\kappa=\kappa(m,n)$ is a positive constant and $\gamma$ is as in~\ref{def-g-+}.  In view of above facts
that $\alpha_{0} r< \varphi(r)\,\forall r$ and $\varphi(r)-\alpha_{0} r\le C r^{\gamma}$ for $r\ge 1$ (hence
$\varphi(r/t)-\alpha_{0} r/t\le C(r/t)^{\gamma}$ for $r\ge t$), we see
that there is $C=C(m,n)$ with
\[%
\varphi_{t}(r)-\alpha_{0} r\le Ct^{1+|\gamma |}(r+t)^{-|\gamma|}\,(\,\le Ct), \,\, \forall r>0,\,\,t>0, %
\dl{vph-lam-bds}
\]%
where $\varphi_{t}(r)=t\varphi(r/t)$.

We shall also need the fact, proved in~\cite[Lemma 7.5]{Sim21a}, that $S=SG(\varphi)$ is strictly stable, in the
sense that there is $\kappa=\kappa(m,n)>0$ such that
\[%
\kappa\int_{S}(\tilde r^{-2}\zeta^{2}+|\nabla_{S}\zeta|^{2})\,d\mu(x,\xi,y)\le %
\int_{S}\bigl(\bigl|\nabla_{S}\zeta\bigr|^{2}-|A_{S}|^{2}\zeta^{2}\bigr)\,d\mu %
\dl{strict-stab-S}
\]%
for $\zeta=\zeta(x,y)\in C_{c}^{1}(\R^{n+m})$, where $\tilde r=|x|+|\xi|$, $|A_{S}|$ is the length of the second fundamental
form of $SG(\varphi)$, and $\mu$ is $(n+m+\ell-1)$-dimensional Hausdorff measure in $\R^{n+m+\ell}$. In fact in
\cite{Sim21a} a weaker inequality with only the term $\tilde r^{-2}\zeta^{2}$ on the left was established, but since
$|A_{S}|^{2}\le C\tilde r^{-2}$ (from \ref{2nd-ff-s}) and $\int_{S}\tilde r^{-2}\zeta^{2}\le C\int_{S}|\nabla\zeta|^{2}$, the
inequality \ref{strict-stab-S} follows directly from this weaker inequality, because, for $\theta\in (0,1)$, the weaker inequality
implies
\[%
\begin{aligned}%
(1-\theta)\kappa\int_{S}\tilde r^{-2}\zeta^{2}&\le
(1-\theta)\int_{S}\bigl(\bigl|\nabla_{S}\zeta\bigr|^{2}-|A_{S}|^{2}\zeta^{2}\bigr)\,d\mu  \\  %
\noalign{\vskip-3pt}
&\le \int_{S}\bigl(\bigl|\nabla_{S}\zeta\bigr|^{2}-|A_{S}|^{2}\zeta^{2}\bigr)\,d\mu 
 -\theta\int_{S}\bigl|\nabla_{S}\zeta\bigr|^{2} +C\theta\int_{S}\tilde r^{-2}\zeta^{2}.
\end{aligned}%
\]%

\vskip3pt

{\bf{}\tl{rem-on-m'} Remark:} If $\tilde m,\tilde n$ (fractional) are sufficiently close to $m,n$ respectively, the
above arguments, including the calibration argument of~\cite{Law72}, apply equally well if we consider the
modified area functional
\[%
\widetilde{\!\!\mathcal{F}}(u)={\textstyle\int_{0}^{1}}\sqrt{1+(u'(r))^{2}}\,u^{\tilde m-1}r^{\tilde n-1}\,dr
\]%
in place of the original~\ref{m-n-fnl}; the Euler-Lagrange equation for this modified functional is the ODE
\[%
\bigl(1+(\varphi')^{2}\bigr)^{-1}\varphi'' +\bigl((\tilde n-1)/r\bigr)\varphi' = (\tilde m-1)/\varphi. %
\leqno{(\ddag)}
\]%
Thus, with $\tilde m,\tilde n$ sufficiently close to $m,n$ respectively, there is a unique solution subject to the
initial conditions $\varphi(0)=1$ and $\varphi'(0)=0$, and this solution satisfies all of the
conditions~\ref{props-varphi} and~\ref{vph-lam-bds} with $\tilde m,\tilde n$ in place of $m,n$, with
$\smash{\bigl(\frac{\tilde m-1}{\tilde n-1}\bigr)^{1/2}}$ in place of $\alpha_{0}$, and with $\tilde \gamma$ in place
of $\gamma$, where
\[%
\tilde\gamma=-\frac{\tilde m+\tilde n-3}{2}+\bigl(\bigl(\frac{\tilde m+\tilde n-3}{2}\bigr)^{2}-  %
(\tilde m+\tilde n-2)\bigr)^{1/2}.  %
\]%

\section{Families of SME Supersolutions} \label{super-solns}

Let $\eta>0$ and
\[%
\tilde{n}-1=(n-1)/(1+\eta),\,\,\tilde{m}-1=(m-1)/(1+\eta).
\]%
We assume for the remainder of the discussion that
\[%
\eta=\eta(m,n)\in (0,\tfrac{1}{4}]  %
\dl{def-eta}
\]%
sufficiently small to ensure that we can select $\tilde\varphi\in C^{\infty}[0,\infty)$, in accordance with the
discussion of Remark~\ref{rem-on-m'}, to satisfy~\ref{rem-on-m'}\hskip1pt$(\ddag)$ and all the
conditions~\ref{props-varphi}, with $\alpha_{0}=\sqrt{\frac{\tilde{m}-1}{\tilde{n}-1}}=\sqrt{\frac{m-1}{n-1}}$, with
$\tilde\varphi$ in place of $\varphi$, and with $\tilde\gamma$ in place of $\gamma$, where
\[%
\tilde\gamma=-\frac{\tilde{m}+\tilde{n}-3}{2}+ %
\sqrt{\bigl(\frac{\tilde{m}+\tilde{n}-3}{2}\bigr)^{2} - (\tilde{m}+\tilde{n}-2)\,\,}. %
\dl{gamma'}
\]%
Notice that then, since $-\frac{x}{2}+\sqrt{\frac{x^{2}}{4}-x-1}$ is increasing in the variable $x$,
\[%
\tilde\gamma< -\frac{m+n-3}{2}+\smash{\bigl(\bigl(\frac{m+n-3}{2}\bigr)^{2} - (m+n-2)\bigr)^{1/2}}=\gamma,
\]%
so the solution $\tilde\varphi(r)$ of Remark~\ref{rem-on-m'} decays to $\alpha_{0} r$ faster than the solution
of~\ref{varphi-equn} as $r\to \infty$. 

Also $(1+\eta)^{-1}\bigl((1+\eta)\frac{\tilde\varphi''}{1+(\tilde\varphi')^{2}} +\bigl(\frac{(n-1)}{r}\tilde\varphi' -
\frac{(m-1)}{\tilde\varphi}\bigr)\bigr) =\frac{\tilde\varphi''}{1+(\tilde\varphi')^{2}}
+\frac{(\tilde{n}-1)}{r}\tilde\varphi' - \frac{(\tilde m-1)}{\tilde\varphi}=0$, so
\[%
\frac{(n-1)}{r}\tilde\varphi' - \frac{(m-1)}{\tilde\varphi} = %
-(1+\eta)\frac{\tilde\varphi''}{1+(\tilde\varphi')^{2}}, %
\dl{vph-sup}
\]%
so in particular $\tilde\varphi$ is a supersolution of the SME on $\R^{n}$.  

For each $\epsilon\in (0,\tfrac{1}{2}]$ we can choose $\tilde\epsilon\in (0,\epsilon^{4})$ such that
$\varphi_{\tilde\epsilon}(1) \le\tilde\varphi_{\epsilon^{4}}(1)$ . Then $\{\tilde\varphi_{t}(r)\}_{t\ge \epsilon^{4}}$ is a family
of supersolutions which are $\ge\varphi_{\tilde\epsilon}(r)$ at $r=1$, and so if $\varphi_{\tilde\epsilon}(r)>
\tilde\varphi_{\epsilon^{4}}(r)$ for some $r<1$ we could select $t>\epsilon^{4}$ with
$\tilde\varphi_{t}-\varphi_{\tilde\epsilon}$ having a zero minimum in $r<1$, contradicting the maximum principle discussion
following~\ref{lin-diff}. Thus
\[%
\varphi_{\tilde\epsilon}(r) \le \tilde\varphi_{\epsilon^{4}}(r) \text{ for all }r\le 1,\,\, \epsilon\in(0,\tfrac{1}{2}]. %
\dl{comp}
\]%
In the following lemma we use the particular supersolution $\tilde\varphi$ to prove the existence of a large family of
supersolutions of the SME on suitable domains in $\R^{n+\ell}$. Here $h$ is as in~\ref{def-h} and, for $\epsilon\ge 0$ and
$t\ge 0$, we let
\[%
\psi_{t,\epsilon}(y)= t+  e^{-1/(\epsilon +h(y))}
\dl{def-h-e}
\]%
and
\[%
\Omega_{\epsilon}=\bigl\{(x,y)\in \R^{n}\times\R^{\ell}:|x|<(\epsilon+h(y))^{2}\bigr\}. %
\dl{def-om-e}
\]%
Note that then
\[%
\Omega_{0}=\Omega\,\,\, \text{($\Omega$ as defined in~\ref{def-om}).}   %
\dl{om-0}
\]%

\begin{state}{\bf{}\tl{sup-soln} Lemma (A Family of Supersolutions.)}%
There is $\tau_{0}=\tau_{0}(\ell,m,n)\in (0,\tfrac{1}{2}]$ such that if $\epsilon,\tau\in (0,\tau_{0}]$,
if~\emph{\ref{def-h}} holds,  and if
\[%
S_{t,\epsilon}(x,y)=\psi_{t,\epsilon}(y)\tilde\varphi(|x|/\psi_{t,\epsilon}(y)), %
\]%
where $\psi_{t,\epsilon}$ is as in~\emph{\ref{def-h-e}} above, then
\[%
\mathcal{M}(S_{t,\epsilon}) <0 \text{ on } \wbar\Omega_{\epsilon},\,\,\,\forall\,t\ge 0.
\]%
\end{state}%

{\bf{}\tl{sup-soln-bd} Remark.}  Note that by \ref{vph-lam-bds} and definition~\ref{def-h-e},
\[%
0<\tilde\varphi_{t}(|x|)-\alpha_{0}|x|\le S_{t,\epsilon}(x,y)-\alpha_{0} |x| %
\le C\psi_{t,\epsilon}(y)\le  C t+C_{\!j}(\epsilon+ h)^{j}
\]%
$\forall\,(x,y)\in\R^{n}\times\R^{\ell},\,j\ge 1,\,t>0$, where $C=C(\ell,m,n)$ and $C_{j}=C(\ell,m,n,j)$, because
$t<\psi_{t,\epsilon}(y)\le t+C_{\!j}\,(\epsilon+ h)^{j}$ for each $j=1,2,\ldots$ and each $t\ge 0$.

\smallskip

\begin{proof}{\bf{}Proof of~\ref{sup-soln}.}
Let $\psi\in C^{\infty}(\R^{\ell})$ with $0<\psi\le 1$ and let
\[%
S(x,y)=s(r,y) = \psi(y)\tilde\varphi\bigl(r/\psi(y)\bigr),\,\,\, r=|x|,\,\,y\in \R^{\ell}.  %
\pdl{def-s}
\]%
Then, with $\mathcal{M}$ as in~\ref{Mu}, $s_{y}=D_{\!y}s=(D_{\!y_{1}}s,\ldots,D_{\!y_{\ell}}s)$, and
$s_{yy}=(s_{y_{i}y_{j}})=(D_{\!y_{i}}D_{\!y_{j}}s)$,
\begin{align*}{%
\mathcal{M}(S) &=  s_{rr} + \tfrac{n-1}{r}s_{r}-\tfrac{m-1}{s} +\Delta_{y}s - %
\tfrac{s_{r}^{2}s_{rr}+\tsum_{i,j=1}^{\ell}s_{y_{i}}s_{y_{j}}s_{y_{i}y_{j}} %
+2s_{r}\tsum_{j=1}^{\ell}s_{y_{j}}s_{ry_{j}}} {1+s_{r}^{2}+|s_{y}|^{2}} %
\ptg{M0-S} \\  %
\noalign{\vskip2pt} 
&\hskip-.2in=  \tfrac{1+|s_{y}|^{2}}{1+s_{r}^{2}+|s_{y}|^{2}}s_{rr} + \tfrac{n-1}{r}s_{r}-\tfrac{m-1}{s} +
\Delta_{y}s - \tfrac{\tsum_{i,j=1}^{\ell}s_{y_{i}}s_{y_{j}}s_{y_{i}y_{j}}+2s_{r}\tsum_{j=1}^{\ell}s_{y_{j}}s_{ry_{j}}} %
{1+s_{r}^{2}+|s_{y}|^{2}}
}\end{align*}%
Since $s_{r}=\tilde\varphi'(r/\psi)$ and $s_{rr}=\psi(y)^{-1}\tilde\varphi''(r/\psi)$, we have by~\ref{vph-sup}
\[%
\tfrac{(n-1)s_{r}}{r}-\tfrac{m-1}{s}=
\tfrac{1}{\psi(y)}\Bigl(\tfrac{(n-1)\tilde\varphi'(r/\psi)}{r/\psi}-\tfrac{(m-1)}{\tilde\varphi(r/\psi)} \Bigr)
=-\tfrac{(1+\eta) \tilde\varphi''(r/\psi)}{(1+s_{r}^{2})\psi(y)}
\]%
with $\eta=\eta(m,n)>0$, so~\ref{M0-S} gives
\begin{align*}{%
\mathcal{M}(S) &= %
 \Bigl(\tfrac{1+|s_{y}|^{2}}{1+s_{r}^{2}+|s_{y}|^{2}}-\tfrac{1+\eta}{1+s_{r}^{2}}\Bigr)s_{rr}  +\Delta_{y}s - %
\tfrac{\tsum_{i,j=1}^{\ell}s_{y_{i}}s_{y_{j}}s_{y_{i}y_{j}}+2s_{r}\sum_{j=1}^{\ell}s_{y_{j}}s_{ry_{j}}} %
{1+s_{r}^{2}+|s_{y}|^{2}} %
\ptg{M0-S-2} \\  %
\noalign{\vskip2pt}
&\le \tfrac{-\eta  + |s_{y}|^{2}}{1+s_{r}^{2}}s_{rr} +\ell|s_{yy}|  %
+\tfrac{|s_{y}|^{2}|s_{yy}|+2s_{r}|s_{y}||s_{ry}|} %
{1+s_{r}^{2}+|s_{y}|^{2}}
}\end{align*}%
so
\[%
(1+s_{r}^{2})\mathcal{M}(S) \le (-\eta+|s_{y}|^{2})s_{rr} +\bigl((\ell+1)|s_{yy}|+ |s_{ry}|\bigr)(1+s_{r}^{2}). %
\pdl{M-S-3}
\]%
Now
\[%
\begin{aligned}%
&s_{r}=\tilde\varphi'(r/\psi), \,s_{y_{j}} =\Phi(r/\psi)\psi_{y_{j}}, \text{ where $\Phi(t) = %
  \tilde\varphi(t)-t\tilde\varphi'(t)$,} \\ %
&s_{rr}=\psi^{-1}\tilde\varphi''(r/\psi),\, s_{ry_{j}}=-r\psi^{-2}\psi_{y_{j}}\tilde\varphi''(r/\psi),\\ %
&s_{y_{i}y_{j}}=\Phi(r/\psi)\psi_{y_{i}y_{j}} + r^{2}\psi^{-3}\psi_{y_{i}}\psi_{y_{j}}\tilde\varphi''(r/\psi).
\end{aligned}%
\]%
By~\ref{props-varphi} and~\ref{rem-on-m'}, there  are constants $k=k(\C_{0}),\,b=b(\C_{0})$ such that
{\abovedisplayskip3pt\belowdisplayskip0pt%
\[%
0<\Phi(t)\le k \text{ and also }\Phi(t)\le b(1+t^{2})\tilde\varphi''(t)\,\,\forall t\ge 0, 
\]}%
so \ref{M-S-3} gives  
\begin{align*}{%
&\tfrac{(1+s_{r}^{2})\psi(y)}{\tilde\varphi''(r/\psi)}\mathcal{M}(S)  \le -\eta+k^{2}|\psi_{y}|^{2}\\ %
\noalign{\vskip-4pt}  %
&\hskip1.1in  +(1+\alpha_{0}^{2}) (\ell+1)\bigl(\tfrac{r^{2}}{\psi^{2}}|\psi_{y}|^{2}
+b\psi(1+\tfrac{r^{2}}{\psi^{2}})|\psi_{yy}|\bigr) +(1+\alpha_{0}^{2})\tfrac{r}{\psi}|\psi_{y}|,%
}\end{align*}%
and so on $\Omega_{\epsilon}=\{(x,y):|x|<h_{\epsilon}^{2}\}$, where $h_{\epsilon}=\epsilon+h$, we have
\begin{align*}{%
&\tfrac{(1+s_{r}^{2})\psi(y)}{\tilde\varphi''(r/\psi)}\mathcal{M}(S)  \le -\eta+k^{2}|\psi_{y}|^{2} \\ %
\noalign{\vskip-5pt}  %
&\hskip1.0in  +(1+\alpha_{0}^{2}) (\ell+1)\bigl(\tfrac{h_{\epsilon}^{4}}{\psi^{2}}|\psi_{y}|^{2}+ %
b\psi(1+\tfrac{h_{\epsilon}^{4}}{\psi^{2}}) |\psi_{yy}|\bigr)+(1+\alpha_{0}^{2})\tfrac{h_{\epsilon}^{2}}{\psi}|\psi_{y}|)%
}\end{align*}%
Thus $\mathcal{M}(S)<0$ on $\Omega_{\epsilon}$ if $\psi(y)>0$ is chosen so that 
\[%
\bigl(k^{2} +(1+\alpha_{0}^{2})  (\ell+1)\tfrac{h_{\epsilon}^{4}}{\psi^{2}} \bigr)|\psi_{y}|^{2} + (1+\alpha_{0}^{2}) %
(\ell+1)b\psi(1+\tfrac{h_{\epsilon}^{4}}{\psi^{2}})|\psi_{yy}|+(1+\alpha_{0}^{2})\tfrac{h_{\epsilon}^{2}}{\psi}|\psi_{y}| <\eta. %
\pdl{M-S-4}   %
\]%
One can now directly check that, if $\tau_{0}=\tau_{0}(\ell,m,n)$ is small enough, and if we take $\epsilon,\tau, h,
\psi=\psi_{t,\epsilon}$ as in the statement of the lemma, then~\ref{M-S-4} does hold, so $\mathcal{M}(S_{t,\epsilon})<0$ on
$\wbar\Omega_{\epsilon}$ as required.
 \end{proof} 

For later reference observe that, with $S_{\epsilon}=S_{\epsilon^{4},\epsilon}$ (i.e.\ $S_{t,\epsilon}$
with $t=\epsilon^{4}$) we have, by~\ref{vph-lam-bds} and the definition~\ref{def-h-e}, with suitable $C=C(m,n)$,
\begin{align*}{%
 h^{-2}_{\epsilon}(y)S_{\epsilon}(0,y) &\le Ch^{-2}_{\epsilon}(y)\psi_{\epsilon^{4},\epsilon}(y) %
\le Ch^{-2}_{\epsilon}(y)  \bigl(\epsilon^{4}+e^{-1/(\epsilon+h(y))}\bigr)  %
\dtg{r=0-bd} \\  %
&\le Ch^{-2}_{\epsilon}(y)h^{4}_{\epsilon}(y) =
 Ch^{2}_{\epsilon}( y) \le C(\epsilon+\tau)^{2} \,\,\,\forall y\in U. %
}\end{align*}%

\section[Solutions $u$ of the SME with Small \smash{$|D_{\!y}u|$}]{Solutions \emph{u} of the SME
with Small \emph{D}$_{\!\normalsize\bf\emph{y}}$\emph{u}} \label{small-D-y-sec}

In this section we establish some conditions for a good $C^{2}$ approximation of the $y=$ const.\ slices of $u$,
plus stability consequences, in case $u-\alpha_{0}r$ is a solution of the SME satisfying a $|D_{\!y}u|$ smallness condition.

We shall need the following consequence of the Liouville-type result established in \cite[Corollary~1]{Sim21a}:

\begin{state}{\bf{}\tl{liouv-co} Lemma.}%
There is $\delta_{0}=\delta_{0}(\ell,m,n)>0$ such that if $u=u(x,y)\in C^{2}(\R^{n+\ell})$ is a positive solution of
the \emph{SME} with $u(x,y) >\alpha_{0} |x|$ everywhere on $\R^{n+\ell}$ and $\max|D_{\!y}u|\le \delta_{0}$,
then $u(x,y)=\varphi_{\lambda}(x)$ for some $\lambda>0$. (In particular $u(x,y)$ is independent~of~$y$.)
\end{state}%

\begin{proof}{\bf{}Proof:} By scaling we can assume $u(0,0)=1$.  Let $u_{R}(x,y)=R^{-1}u(Rx,Ry)$ for $R>1$. Then
$u_{R}(0,0)=R^{-1}\to 0$ as $R\to\infty$, so by Lemma~\ref{il-sw-app} $u_{R}$ converges locally uniformly to
$\alpha_{0}|x|$ on $\R^{n+\ell}$.  So $M=SG(u)$ has $\C$ (as in~\ref{def-C}), with multiplicity~$1$, as its tangent cone at
$\infty$, and hence $M=SG(u)$ satisfies the hypotheses of Corollary~1 of~\cite{Sim21a} and so $u(x,y)=\varphi(|x|)$.
\end{proof}

In the following theorem we establish $C^{2}$ bounds and strict stability for solutions $u$ with $u-\alpha_{0}|x|>0$.

In the statement of the theorem we let $V$ be any open subset of $\R^{\ell}$, $q$ Lipschitz with $q>0$ on $V$, and we
suppose
\[%
\left\{\begin{aligned}%
&\,\, W=\{(x,y):|x|<q(y),\,y\in V\}  \\  %
&\,\, |q(y_{1})-q(y_{2})|\le \tau|y_{1}-y_{2}|, \quad y_{1},y_{2}\in V %
\end{aligned}\right.%
\dl{bds-q}
\]%
for a given $\tau>0$ (to be chosen).

\begin{state}{\bf{}\tl{C-2-approx} Theorem (\emph{C}$^{\bf 2}$ approximation.)}%
Let $\delta\in (0,\delta_{0}]$ with $\delta_{0}=\delta_{0}(\ell,m,n)$ as in~\emph{\ref{liouv-co}}.  There is
$\tau=\tau(\ell,m,n,\delta)\in (0,\tfrac{1}{2}]$ such that if~\emph{\ref{bds-q}} holds, if $u$ satisfies the \emph{SME} on $W$
with $0<u(x,y)-\alpha_{0}|x|<\tau q(y)$ and $|D_{\!y}u(x,y)|\le\delta_{0}$ on $W$, and
$|D(u(x,y)-\alpha_{0}|x|)|+q(y)\bigl|D^{2}(u(x,y)-\alpha_{0}|x|)\bigr|<\tau$ for $\tfrac{1}{2}q(y)<|x|<q(y)$, then, with
$\lambda_{y}=u(0,y)$ and $\varphi_{\lambda_{y}}(r)=\lambda_{y}\varphi(r/\lambda_{y})$,
\begin{align*}{%
&(\lambda_{y}+|x|)^{-1}|u(x,y)-\varphi_{\lambda_{y}}(|x|)| + |D_{x,z}(u(x,z)-\varphi_{\lambda_{y}}(|x|))\bigl|_{z=y}  %
\tg{\rm (i)}\\ %
\noalign{\vskip-3pt}
&\hskip1.5in + (\lambda_{y}+|x|) |D_{x,z}^{2}(u(x,z)-\varphi_{\lambda_{y}}(|x|))\bigl|_{z=y}<\delta %
}\end{align*}%
for all $(x,y)\in W$, where $D_{x,z}=(D_{x},D_{z})$, $D_{x,z}^{2}=(D_{x}^{2},D_{x}D_{z},D_{z}^{2})$ (in particular
$|D_{\!y}u(x,y)|<\delta$), and for each $R>0$ we have the strict stability inequality
\[%
\kappa\int_{M_{R}}(|\nabla_{M}\zeta|^{2}+\tilde r^{-2}\zeta^{2})\,d\mu(x,\xi,y) %
\le \int_{M_{R}}\bigl( |\nabla_{M}\zeta|^{2} - |A_{M}|^{2}\zeta^{2}\bigr)\,d\mu %
\leqno{\rm (ii)}
\]%
for all $\zeta=\zeta(x,y)\in C_{c}^{1}(W)$, where $M=SG(u)$, $M_{R}= M\cap \{(x,y)\in
\R^{n+\ell}:|y^{j}|<R,\,j=1,\ldots,\ell\}$, and $\kappa=\kappa(\ell,m,n)>0$.
\end{state}%

{\bf{}\tl{rem-on-shape} Remarks} {\bf{}(1)} \,Notice we just need $\zeta=0$ on $\partial W$ in~(ii); there is no necessity that
$\zeta=0$ at the points $|y^{j}|=R$. Also we shall see below that the proof of~(ii) makes no use of the fact that $u$ satisfies
the SME; we just need~(i) with small enough $\delta$ to prove the strict stability~(ii).

{\bf{}(2)} \,Assuming $u_{\tau}$ as in Theorem~\ref{ex-thm} exists (to be proved in \S\hskip1.5pt\ref{soln-family}) the above
applies with $V=U$, $q=\frac{3}{4}h^{2}$ and $u=u_{\tau}$ , so (i) holds in this case, and then $u$ as in~\ref{def-u},
\ref{props-v} also satisfies~(i) in case $q=R$ (constant) for any choice of $R>1$ and, in view of (1) above, we have the strict
stability~(ii) for $SG(u)$ with $u$ as in~\ref{def-u}.  Also, since $f$ is smoothly as close
to~$1$ as we wish, this also gives the required strict stability of $SG(u)$, $u$ as in~\ref{def-u}, with respect to the metric
$\sum dx_{i}^{2}+ f(x,y)\sum d\xi_{j}^{2}+\sum_{}dy_{k}^{2}$.

{\bf{}(3)} \,Notice that the fact that an inequality like (i) holds for $u$ as in~\ref{def-u} means that we have a fairly precise
picture of the shape of $M=SG(u)$: For each $y_{0}\in \R^{\ell}$ the slice $M\cap \{(x,y):y=y_{0}\}$ of $M=SG(u)$ is
$SG(\alpha_{0}r)$ if $y_{0}\in K$ while if $y_{0}\in U\,\,(=\R^{\ell}\setminus K)$ the slice, after rescaling, is $C^{2}$ close
to $SG(\varphi)$.

\smallskip

\begin{proof}{\bf{}Proof of~\ref{C-2-approx} (i):}   
It suffices to prove~(i) with $y=0$ and, by rescaling, we can assume $q(0)=1$. Then, since and $|q(y)-q(0)|\le \tau|y|$, the
domain of $u\supset B_{1/2}$ for $\tau\le \tfrac{1}{2}$. 

We first show that, with $\lambda=u(0,0)$, there is $\tau=\tau(\ell,m,n,\delta)$ such that the stated hypotheses imply
\[%
(\lambda+|x|)^{-1}|u(x,y)-\varphi_{\lambda}(|x|)|<\delta,\quad |x|\le q(y), \,\,|y|\le |x|. 
\pdl{c-2-a-1}
\]%
First observe that in case $\tau\le \tfrac{1}{20}\delta$, \ref{c-2-a-1} holds trivially if $\tfrac{1}{8}q(y)\le |x|\le q(y)$ by virtue of
the given inequality $u(x,y)-\alpha_{0}|x|\le \tau q(y)$ and the fact that $\varphi_{\lambda}(r)-\alpha_{0}r\le \lambda\le \tau$. 
So it suffices to prove~\ref{c-2-a-1} in case $|x|<\tfrac{1}{8}q(y)$, so it is in fact sufficient to prove~\ref{c-2-a-1} for
$|x|\le\tfrac{1}{4}$ (because $|x|>\tfrac{1}{4}\Rightarrow |x|>\tfrac{1}{8}q(y)$ if $|y|\le |x|$ and $\tau\le \tfrac{1}{4}$). To do
this, we first claim that \ref{c-2-a-1} also holds for $\tfrac{\lambda}{\eta}\le |x|\le \tfrac{1}{4}$ if $\tau<\tfrac{1}{4}\eta$ (hence
$\lambda<\tfrac{1}{4}\eta$), where $\eta$ is the constant $\eta(\ell,m,n,\theta,\delta)$ of Lemma~\ref{il-sw-app} with
$\theta=\tfrac{1}{2}$ and with $\tfrac{1}{8}\delta$ in place of $\delta$ (so in particular $\eta\le \tfrac{1}{8}\delta$). In that
case $\tfrac{\lambda}{\eta}<\tfrac{1}{4}$ and we can apply~\ref{il-sw-app} to the scaled function
$u_{\rho}(x,y)=\rho^{-1}u(\rho x,\rho y)$ with $\tfrac{\lambda}{\eta}\le\rho$ (so $u_{\rho}(0,0)=\lambda/\rho\le \eta$ and
the domain of $u_{\rho}$ contains $\breve B_{1}$). Thus~\ref{il-sw-app} gives $u_{\rho}(x,y)-\alpha_{0}|x|\le
\tfrac{1}{8}\delta$ for all $(x,y)\in B_{1/2}$. In terms of the original $u$ this gives
\[%
\rho^{-1}(u(\rho x,\rho y)-\alpha_{0}|\rho x|) \le \tfrac{1}{8}\delta, \quad (x,y)\in B_{1/2}.
\pdl{il-sw-conseq}
\]%
If $|y_{0}|\le |x_{0}|$ and $\tfrac{\lambda}{\eta}\le|x_{0}|\le \tfrac{1}{4}$ then with $\rho=4|x_{0}|$ and
$(x,y)=\tfrac{1}{4}|x_{0}|^{-1}(x_{0},y_{0})\,\,(\in B_{1/2})$,~\ref{il-sw-conseq} implies
$|x_{0}|^{-1}(u(x_{0},y_{0})-\alpha_{0}|x_{0}| )\le \tfrac{1}{2}\delta$, and since $|\alpha_{0}r-\varphi_{\lambda}(r)|\le
\lambda$ this gives $|x_{0}|^{-1}|u(x_{0},y_{0})-\varphi_{\lambda}(|x_{0}|)| \le \tfrac{1}{2}\delta + \tfrac{\lambda}{|x_{0}|}<
\delta$ because $|x_{0}|\ge \lambda/\eta$.  So~\ref{c-2-a-1} holds at the point $(x_{0},y_{0})$ and we have checked
that~\ref{c-2-a-1} holds with suitable $\tau=\tau(\ell,m,n,\delta)$ for $\lambda/\eta\le |x|\le q(y)$.

So to complete the proof of~\ref{c-2-a-1} we can assume that $(x,y)\in W$ with $|y|\le |x|\le p\lambda\le \tfrac{1}{2}q(y)$,
where $p\,\,(=1/\eta)$ is fixed.  If~\ref{c-2-a-1} fails in this case then there are sequences $q_{k},u_{k},\tau_{k}$ with
$\tau_{k}\to 0$, \ref{bds-q} holds with $q=q_{k}$, $q_{k}(0)=1$, and $u_{k}$ are SME solutions in
$W_{k}=\{(x,y):|x|<q_{k}(y)\}$ such that the hypotheses of the theorem hold with $\tau=\tau_{k}$, $q=q_{k}$ yet there are
$(x_{k},y_{k})\in W_{k}$ with $|y_{k}|\le|x_{k}|$, $|x_{k}|\le p\lambda_{k}$ and
\[%
(\lambda_{k}+|x_{k}|)^{-1}|u_{k}(x_{k},y_{k})-\varphi_{\lambda_{k}}(|x_{k}|)|\ge \delta,  %
\quad   \lambda_{k}=u_{k}(0,0)\le \tau_{k}\to 0. %
\pdl{c-2-a-2}
\]%
We let $\tilde q_{k}(y)=\lambda_{k}^{-1}q_{k}(\lambda_{k}y)$ and $\tilde u_{k}(x,y)=\lambda_{k}^{-1}u_{k}(\lambda_{k}
x,\lambda_{k} y)$ for $(x,y)\in B_{\lambda_{k}^{-1}/2}$ $(\,\subset \{(x,y):|x|<\tilde q_{k}(y)\}$ for sufficiently large $k$
since $q_{k}(0)=1$ and $|Dq_{k}|\le\tau_{k}$).  For $R>1$, define $\tilde u_{k,R}=R^{-1}\tilde u_{k}(Rx,Ry)$. Then $\tilde
u_{k,R}(0,0)=1/R$ and the domain of $\tilde u_{k,R}\supset\breve B_{1}$. Thus for any $R>\eta_{1}^{-1}$, with $\eta_{1}$
as in~\ref{grad-bd}, and for any $k$ sufficiently large (depending on $R$) we can apply~\ref{il-sw-app} and \ref{grad-bd} to
give a positive $\kappa_{R},C_{R}$ with $\kappa_{R}\le \tilde u_{k,R}-\alpha_{0}|x|\le 1$ and $|D\tilde u_{k,R}|\le C_{R}$
on $B_{1/2}$, so
\[%
R\kappa_{R}\le \tilde u_{k}(x,y)-\alpha_{0}|x|\le R,\,\, |D\tilde u_{k}|  \le C_{R}, \quad R>1,\,\, (x,y)\in B_{R/2},
\]%
for all sufficiently large $k$, where $C_{R},\kappa_{R}$ do not depend on $k$.  Thus $|D\tilde u_{k}|$ is locally uniformly
bounded on $\R^{n+\ell}$ and so a subsequence of $\tilde u_{k}$ converges locally in $C^{2}$ by standard estimates for
solutions of quasilinear elliptic PDE with bounded gradient (see e.g.\ \cite{GilT83}) to a positive solution $u$ of the SME
with $u(0,0)=1$, $u-\alpha_{0}|x|\ge 0$ (hence $>0$ by the maximum principle), and $|D_{y}u|\le \delta_{0}$ everywhere on
$\R^{n+\ell}$.  Then, by Lemma~\ref{liouv-co}, $u(x,y)=\varphi(|x|)$, so in fact $\tilde u_{k}$ converges locally in the
$C^{2}$ sense to $\varphi(|x|)$.  Thus $\tilde u_{k}(x,y)-\varphi(|x|)\to 0$ in $B_{R_{0}}$ uniformly for each $R_{0}>0$,
which in terms of the original $u_{k}$ gives $\lambda_{k}^{-1}u_{k}(\lambda_{k}x,\lambda_{k}y)-\varphi(|x|)=
\lambda_{k}^{-1}(u_{k}(\lambda_{k}x,\lambda_{k}y)-\varphi_{\lambda_{k}}(|\lambda_{k}x|))\to 0$ uniformly on
$B_{R_{0}}$.  $|(x_{k},y_{k})|\le 2p\lambda_{k}$, so we can take $(x,y)=\lambda_{k}^{-1}(x_{k},y_{k})$, and hence
$\lambda_{k}^{-1}|u_{k}(x_{k},y_{k})-\varphi_{\lambda_{k}}(|x_{k}|)| \to 0$, contradicting~\ref{c-2-a-2}.

Thus~\ref{c-2-a-1} is proved, and, in combination with~\ref{grad-est}, this enables us to apply interior quasilinear elliptic
estimates \cite{GilT83} to $u$ and also, by~\ref{lin-eqn-for-diff} (with $u_{1}=u$ and $u_{2}=\varphi_{\lambda}(|x|)$), to the
difference $u-\varphi_{\lambda}$, in the conical domain $|x|< \frac{3}{4}q(y),\, |y|<|x|$, giving~$(\ddag)$ for
$|x|\le\tfrac{1}{2}q(y)$ with $y=0$.  Since 
$|D(u(x,y)-\alpha_{0}|x|)|+q(y)\bigl|D^{2}(u(x,y)-\alpha_{0}|x|)\bigr|<\tau$ for $\tfrac{1}{2}q(y)<|x|<q(y)$ is given, and,
by~\ref{props-varphi}, $\varphi_{\lambda}(r)-\alpha_{0}r=\lambda\Psi(\lambda/r)$, where $|\Psi'(r)|+r|\Psi''(r)|\le C(m,n)$
for $r\ge 1$, we also have the required inequality $(\ddag)$ for $y=0$ in $\tfrac{1}{2}q(y)<|x|<q(y)$.

\setcounter{pequation}{0}

{\bf{}Proof of (ii):} For any $z\in V$ we let $\lambda_{z}=u(0,z)$, $u_{z}(x,y)=u(x,y+z)$ (so $\lambda_{z}=u_{z}(0,0)$), and
$\tilde u_{z}(x,y)=\lambda_{z}^{-1}u((0,z)+\lambda_{z}(x,y))\,\,(\,=\lambda_{z}^{-1}u_{z}(\lambda_{z}(x,y)))$. Then by~(i)
\begin{align*}{%
&(1+|x|)^{-1}|\tilde u_{z}(x,y)-\varphi(|x|)|+ |D_{x,y}(\tilde u_{z}(x,y)-\varphi(|x|))| 
\ptg{ss-2}\\  %
\noalign{\vskip-3pt}
&\hskip1.7in +(1+ |x|) |D_{x,y}^{2}(\tilde u_{z}(x,y)-\varphi(|x|))|<\delta\,\,\, \text{ at $y=z$}. 
}\end{align*}%
Let $\zeta=\zeta(x,y)\in C_{c}^{1}(W)$, and for $z\in V$ let $\zeta_{z}(x)=\zeta(x,z)$
and $\tilde\zeta_{z}(x)=\zeta(\lambda_{z}^{-1}x,z)$.  Also, for $z\in V$ let
\[%
\left\{\begin{aligned}%
M_{z}&=M\cap\{(x,\xi,y):y=z\} \,\,\bigl(= SG(u_{z}\bigl|\{(x,y):y=z\})\bigr)  \\
\wtilde M_{z}&=\lambda_{z}^{-1}(M_{z}-(0,0,z))\,\,\bigl(\,=SG( \tilde u_{z}\bigl|\{(x,y):y=0\})\bigr),
\end{aligned}\right.%
\pdl{def-m-z-til}
\]%
Then~\ref{ss-2} implies
{\abovedisplayskip3pt\belowdisplayskip3pt%
\[%
\int_{\wtilde M_{z}}\bigl|\nabla_{\wtilde M_{z}}\tilde\zeta_{z}\bigr|^{2}  %
=\int_{SG(\varphi)}  \bigl|\nabla_{SG(\varphi)}\tilde\zeta_{z}\bigr|^{2}              +E_{1}, %
\pdl{comp-1} 
\]%
\[%
\int_{\wtilde M_{z}}|A_{\tilde M_{z}}|^{2}\tilde\zeta_{z}^{2}\,d\mu         %
= \int_{SG(\varphi)} |A_{SG(\varphi)}|^{2}\tilde\zeta_{z}^{2}            +E_{2},  %
\pdl{comp-2}  %
\]%
and
\[%
\int_{\wtilde M_{z}}(1+|x|)^{-2}\tilde\zeta_{z}^{2}\,d\mu         %
= \int_{SG(\varphi)} (1+|x|)^{-2}\tilde\zeta_{z}^{2}            +E_{3}, %
\pdl{comp-3} %
\]}%
where 
\[%
|E_{j}| \le C\delta\int_{\wtilde M_{z}}\bigl((1+|x|)^{-2}\tilde\zeta_{z}^{2}+ %
                                 \bigl|\nabla_{\wtilde M_{z}}\tilde\zeta_{z}\bigr|^{2}\bigr)\,d\mu, \quad j=1,2,3.
\]%
Notice that here we used $|D_{x}\tilde\zeta_{z}|^{2}\le C|\nabla_{\wtilde M_{z}}\tilde \zeta_{z}|^{2}$, $C=C(m,n)$, which
holds because by~\ref{ss-2} the gradient of $u$ is bounded by the fixed constant $\alpha_{0}+\delta\le \alpha_{0}+1$.

By~\ref{strict-stab-S},
\[%
\kappa\int_{SG(\varphi)}(\tilde r^{-2}\tilde\zeta_{z}^{2}+|\nabla_{SG(\varphi)}\tilde\zeta_{z}|^{2})\,d\mu \le %
\int_{SG(\varphi)}\bigl(\bigl|\nabla_{SG(\varphi)}\tilde\zeta_{z}\bigr|^{2}-  %
|A_{SG(\varphi)}|^{2}\tilde\zeta_{z}^{2}\bigr)\,d\mu,  %
\pdl{ss-3}
\]%
$\tilde r=|x|+|\xi|$, and using~\ref{comp-1}, \ref{comp-2} and~\ref{comp-3} this yields
\[%
 (\kappa-C\delta)\int_{\tilde M_{z}}(\tilde r^{-2}\tilde\zeta_{z}^{2}+|\nabla_{\tilde M_{z}}\tilde\zeta_{z}|^{2})\,d\mu \le %
\int_{\tilde M_{z}}\bigl(\bigl|\nabla_{\tilde M_{z}}\tilde\zeta_{z}\bigr|^{2}-  %
|A_{\tilde M}|^{2}\tilde\zeta_{z}^{2}\bigr)\,d\mu 
\pdl{ss-4}
\]%
and, after changing back to the original scale,
\[%
 (\kappa-C\delta)\int_{M_{z}}(\tilde r^{-2}\zeta_{z}^{2}+|\nabla_{M_{z}}\zeta_{z}|^{2})\,d\mu \le %
\int_{M_{z}}\bigl(\bigl|\nabla_{M_{z}}\zeta_{z}\bigr|^{2}-  %
|A_{M_{z}}|^{2}\zeta_{z}^{2}\bigr)\,d\mu, 
\pdl{ss-5}
\]%
where $\zeta_{z}(x)=\zeta(x,z)$ and $M_{z}=M\cap \{(x,\xi,y):y=z\}$.

Again using~\ref{ss-2}, we have $|D_{\!y}u|<\delta$, $|A_{M_{z}}-A_{M}|\le C\delta/\tilde r$ on $M_{z}$ and $|A_{M}|\le
C/\tilde r$ ($\tilde r=|x|+|\xi|$),  where $C=C(\ell,m,n)$. So integrating with respect to $z$ in~\ref{ss-5} over
the region $\{(x,z)\in W:|z^{j}|<R,\,\,j=1,\ldots,\ell\}$, and using the coarea formula
together with the fact that $|\nabla_{M_{z}}\zeta_{z}|\le |\nabla_{M}\zeta|$, we conclude
\[%
(\kappa-C\delta)\int_{M_{R}}\tilde r^{-2}\zeta^{2}\,d\mu\le %
\int_{M_{R}}\bigl(\bigl|\nabla_{M}\zeta\bigr|^{2}-|A_{M}|^{2}\zeta^{2}\bigr)\,d\mu
\pdl{ss-7}
\]%
with $C=C(\ell,m,n)$. Since $|A_{M}|^{2}\le C\tilde r^{-2}$ by~(i) this gives the required inequality~\ref{stab-ineq},
with $\tilde \kappa=\tilde \kappa(\ell,m,n)$ in place of $\kappa$, for $\delta$ sufficiently small (depending only on
$\ell,m,n$). \end{proof}

\section{Proof of Theorem \ref{ex-thm}} \label{soln-family}

Let $h$ be as in~\ref{def-h} with $\tau\in(0,\tau_{0}]$, $\tau_{0}$ to be chosen, and let
$\Omega_{\epsilon},S_{t,\epsilon}$ be as in~\S\ref{super-solns}.

\begin{state}{\bf{}\tl{periodic-ex} Lemma.}%
Let $\delta\in (0,\delta_{0})$ with $\delta_{0}=\delta_{0}(\ell,m,n)$ as in~\emph{\ref{liouv-co}}, and (as
in~\emph{\ref{r=0-bd}}) let $S_{\epsilon}=S_{\epsilon^{4},\epsilon}$ (i.e.\ $S_{t,\epsilon}$ as in~\emph{\ref{sup-soln}} with
$t=\epsilon^{4}$). There is $\tau_{0}=\tau_{0}(\ell,m,n,\delta)\in (0,\tfrac{1}{2}]$ such that if $h$ is as in~\emph{\ref{def-h}}
and if $\epsilon,\tau\in(0,\tau_{0}]$, then there is a positive solution $u_{\epsilon}\in
C^{2,\alpha}(\wbar{\Omega}_{\epsilon})$ of the \emph{SME} with $u_{\epsilon} = S_{\epsilon}$ on
$\partial\Omega_{\epsilon}$ and with
\[%
\begin{cases}%
0<u_{\epsilon}(x,y)-\alpha_{0}|x| \le S_{\epsilon} -\alpha_{0}|x| \,\, %
\le C(\epsilon^{4}+  h_{\epsilon}^{j}(y)), &(x,y)\in\Omega_{\epsilon} \\ %
|Du_{\epsilon}|\le 2\alpha_{0},\,\,\,|D_{\hskip-0.5pt y}u_{\epsilon}|\le \delta, & (x,y)\in\Omega_{\epsilon},   %
\end{cases}%
\leqno{(\ddag)}
\]%
where $C=C(\ell,m,n,j)$, $j=1,2,\ldots$ and $h_{\epsilon}=\epsilon+h$.
\end{state}%

\begin{proof}{\bf{}Proof of~\ref{periodic-ex}:} 
Assume for the moment that $U$ is bounded, let $\Z=\{0,\pm 1,\pm 2,\ldots\}$,
$\Z^{\ell}=\{(z_{1},\ldots,z_{\ell}):z_{j}\in\Z,\,j=1,\ldots,\ell\}$, and for $R>0$ with $\wbar
U\subset\{y:|y^{j}|<R,\,j=1,\ldots,\ell\}$, extend $h\bigl|\{y:|y^{j}|\le R,\,j=1,\ldots,\ell\}$ to an $R$-periodic function
$h_{R}$ on $\R^{\ell}$; thus
\[%
h_{R}(y+R z) = h(y), \quad z\in\Z^{\ell},\,\,  |y^{j}|<R,\,j=1,\ldots,\ell,
\pdl{small-D-y-1d}
\]%
and the conditions \ref{def-h}  still hold with $h_{R}$ in place of $h$.

For $\epsilon\in[0,1)$, we let
\[%
\Omega_{\epsilon,R}=\{(x,y)\in\R^{n}\times\R^{\ell}:|x|<(\epsilon+h_{R}(y))^{2}\}. 
\]%
Then $\Omega_{0,R}$ is the periodic extension $\cup_{z\in\Z^{\ell}}\bigl(Rz+\Omega\bigr)$ of
$\Omega=\{(x,y):|x|<h^{2}(y)\}$, and $\Omega_{\epsilon,R}$ has smooth boundary for each $\epsilon>0$. For
$\epsilon>0$ we let $S_{\epsilon,R}$ be defined as for $S_{\epsilon}$ with $h_{R}$ in place of $h$; thus
\[%
S_{\epsilon,R}(y) = %
(\epsilon^{4}+e^{-1/(\epsilon+h_{R}(y))})\,\tilde\varphi\bigl(|x|/(\epsilon^{4}+e^{-1/(\epsilon+h_{R}(y))})\bigr), %
\quad y\in\wbar\Omega_{\epsilon,R}. %
\]%
We consider the operator $\mathcal{H}(u)=\sum_{i,j=1}^{n}\bigl(\delta_{ij}-\nu_{i}(u)\nu_{j}(u)\bigr)D_{i}D_{j}u$ as
in~\ref{mco}. Then the SME $\mathcal{M}(u)=0$ can be written,  as in~\ref{Mu}, 
\[%
\mathcal{H}(u) = \tfrac{m-1}{u}. 
\]%
For $\alpha\in(0,1)$ given, let $C^{2,\alpha}_{R}(\wbar{\Omega}_{\epsilon,R})$ denote the $R$-periodic
$C^{2,\alpha}$ functions on $\wbar\Omega_{\epsilon,R}$:
\[%
C^{2,\alpha}_{R}(\wbar{\Omega}_{\epsilon,R})=\{u\in C^{2,\alpha}(\wbar{\Omega}_{\epsilon,R}): %
u(x,y+Rz)=u(x,y)\,\forall\,z\in\Z^{\ell},\,(x,y)\in\wbar{\Omega}_{\epsilon,R}\}.  %
\pdl{c-2-a-q} %
\]%

For $\epsilon>0$ and $\sigma\in[0,1]$, we consider the Dirichlet problem of finding $u=u_{\sigma}\in
C^{2,\alpha}_{R}(\wbar\Omega_{\epsilon,R})$ with
\[%
\begin{cases}%
\mathcal{H}(u) = \tfrac{m-1}{u} &\text{ in }\Omega_{\epsilon,R} \\
u = \varphi_{\tilde\epsilon}+\sigma(S_{\epsilon,R}-\varphi_{\tilde\epsilon}) &\text{
on }\partial\Omega_{\epsilon,R}, 
\end{cases}%
\pdl{dir-prob}
\]%

Notice that the function $u=\varphi_{\tilde\epsilon}\,\,(\,\le S_{\epsilon,R})$ is a suitable solution in case $\sigma=0$.

Suppose that $\sigma_{0}\in [0,1]$ is any value of $\sigma$ such that \ref{dir-prob} has a solution with $|Du_{\sigma_{0}}|\le
2\alpha_{0}$ and $|D_{y}u|\le \delta_{0}$ with $\delta_{0}=\delta_{0}(\ell,m,n)$ as in~\ref{C-2-approx}.  For any such
$\sigma_{0}\in (0,1]$ we automatically have the bounds
\[%
\varphi_{\tilde\epsilon} \le u_{\sigma_{0}}(x,y) \le S_{\epsilon,R} %
\text{ on }\Omega_{\epsilon,R}. %
\pdl{bds-u}
\]%
To check this note that if $u_{\sigma_{0}}<\varphi_{\tilde\epsilon}$ at some point of
$\Omega_{\epsilon,R}$ then we can select $t<\tilde\epsilon$ such that $\varphi_{t}-u_{\sigma_{0}}$ has a zero maximum
$\Omega_{\epsilon,R}$, and this contradicts the maximum principle because, using~\ref{lin-eqn-for-diff},
$\mathcal{L}_{\varphi_{t},u_{\sigma_{0}}}(\varphi_{t}-u_{\sigma_{0}})\ge 0$.  Similarly if $S_{\epsilon,R}<u_{\sigma_{0}}$ at
some point of $\Omega_{\epsilon}$, then we get a contradiction by using the minimum principle to
$S_{t,\epsilon,R}-u_{\sigma_{0}}$ ($S_{t,\epsilon,R}=S_{t,\epsilon}$ as in~\ref{sup-soln} with $h_{R}$ in place of $h$) with
$t\ge \epsilon^{4}$ chosen so that $S_{t,\epsilon,R}-u_{\sigma_{0}}$ has a zero minimum in $\Omega_{\epsilon}$.

We also claim that, for any such $u_{\sigma_{0}}$ with $|D_{y}u_{\sigma_{0}}|\le\delta_{0}$, the inequality
\ref{C-2-approx}\,(i) and also the strict stability inequality~\ref{C-2-approx}\,(ii) holds with $V=U$ and
$q=h_{\epsilon,R}^{2}$.  To check this let $y_{0}\in\R^{\ell}$ be arbitrary and note that, by~\ref{bds-u} and \ref{r=0-bd},
\[%
h_{\epsilon,R}(y_{0})^{-1}u_{\wbar{\sigma}}(0,y_{0})\le C(\epsilon+\tau_{0})^{2},\quad C=C(\ell,m,n).
\]%
So with $\epsilon\le \tau_{0}$, and $\tau_{0}=\tau_{0}(n,m,\ell,\delta)$ small enough, we can
apply~\ref{C-2-approx} as claimed.  So in particular, in the notation of~\ref{jac-0}, the first eigenvalue of the
linearized operator $-\mathcal{L}_{u_{0}}(v)=-V^{-1}\frac{d}{dt}\bigl|_{t=0}\bigl(\mathcal{M}(u_{\sigma_{0}}+tv)\bigr)$ (as an
operator on  $v\in C_{R}^{2,\alpha}(\wbar\Omega_{\epsilon,R})$ with $v=0$ on
$\partial\Omega_{\epsilon,R}$) is $>0$, hence $\mathcal{L}_{u_{0}}$ has no  eigenvalue $=0$. 

So if $\sigma_{0}<1$, we can apply the contraction mapping principle in a neighborhood of $u_{\sigma_{0}}$ to prove
that~\ref{dir-prob} has a solution $u_{\sigma}\in C^{2,\alpha}(\wbar\Omega_{\epsilon,R})$ with $|Du_{\sigma}|\le\alpha_{0}
+\delta$ and $|D_{y}u_{\sigma}|\le\delta$ for $\sigma\in [\sigma_{0},\sigma_{0}+\beta]$ for sufficiently small $\beta>0$. 
We emphasise that this is valid for any $\sigma_{0}\in [0,1)$ such that a positive solution $u=u_{\sigma_{0}}$ of
\ref{dir-prob} exists and satisfies $|Du_{\sigma_{0}}|\le 2\alpha_{0}$ and $|D_{y}u_{\sigma_{0}}|\le \delta_{0}$.

Now let
\begin{align*}{%
&\quad\wbar{\sigma}=\sup\bigl\{t\in [0,1]: u_{\sigma}\in  %
                                    C^{2,\alpha}_{R}(\wbar{\Omega}_{\epsilon}) \text{ satisfies~\ref{dir-prob}, } %
\ptg{main-pf-2} \\ %
\noalign{\vskip-3pt}
&\hskip0.3in \text{  and has the properties } \sup|Du_{\sigma}|<2\alpha_{0}, \,\,
\,|D_{y}u_{\sigma}|<\delta_{0}\,\,  \forall \sigma\in[0, t)\bigr\}. %
}\end{align*}
$\wbar\sigma$ is well defined because the set on the right is non-empty by virtue of the fact that we can apply the
above discussion with $\sigma_{0}=0$ and $u_{0}=\varphi_{\tilde\epsilon}$. 

Take any sequence $\sigma_{k}\uparrow \wbar{\sigma}$.  By~\ref{bds-u}, the estimate~\ref{grad-est} is applicable and in
combination with standard quasilinear elliptic estimates~\cite{GilT83} gives a fixed bound on the $C^{2,\alpha}$ norm of
$u_{\sigma_{k}}$, independent of $k$ (but depending on $\epsilon,\tau$).  So a subsequence of $u_{\sigma_{k}}$ converges
in $C^{2}$ to a positive solution $u_{\wbar{\sigma}}\ge \varphi_{\tilde\epsilon}$ of the SME satisfying~\ref{dir-prob} with
$\sigma=\wbar{\sigma}$ and also satisfying~\ref{bds-u} and
\[%
\max|D_{\!y}u_{\wbar{\sigma}}|\le \delta_{0},\,\,\,\max|Du_{\wbar{\sigma}}|\le 2\alpha_{0}. %
\pdl{max-D-y}
\]%
But, applying the discussion above with $\sigma_{0}=\wbar{\sigma}$, by Theorem~\ref{C-2-approx}\,(i) we have
\[%
|Du_{\wbar{\sigma}}|<\alpha_{0}+\delta\text{ and }|D_{\!y}u_{\wbar{\sigma}}(x,y)|< \delta, %
\quad (x,y)\in\wbar\Omega_{\epsilon,R}, %
\pdl{grad-bd-y1}
\]%
and in particular strict inequality holds in both the inequalities in~\ref{max-D-y}.

But then if $\wbar{\sigma}<1$ we can apply the above discussion for $u_{\sigma_{0}}$ with $\sigma_{0}=\wbar{\sigma}$ to
contradict the definition of $\wbar{\sigma}$ in~\ref{main-pf-2}.  So $\wbar{\sigma}=1$ and, by~\ref{max-D-y}
and~\ref{bds-u}, $u=u_{1}$ satisfies~\ref{dir-prob} with $\sigma=1$ together with the bounds~$(\ddag)$.

In the general case when $U$ is allowed to be unbounded, we take a fixed $C^{\infty}$ function $\zeta:\R^{\ell}\to [0,1]$
with $\zeta(y)=1$ for $|y|\le \tfrac{1}{2}$ and $\zeta(y)=0$ for $|y|\ge \tfrac{3}{4}$, let $\zeta_{R}(y)=\zeta(R^{-1}y)$, and
replace $h$ by $\zeta_{R}\,h$.  Then for suitable choice of $\tau_{0}=\tau_{0}(\ell,m,n,\delta)$ (independent of $R$) the
above discussion applies with $\zeta_{R}\,h$ in place of $h$, giving solutions $u_{R}$ on $\Omega_{\epsilon,R}$ with
\[%
\varphi_{\tilde\epsilon}\le u_{R}\le S_{\epsilon,R},\,\, |Du_{R}|\le 2\alpha_{0}\text{ and }|D_{y}u_{R}|\le \delta.  
\pdl{bds-u-R}
\]%
We then let $R\uparrow\infty$ and observe that, by virtue of~\ref{bds-u-R}, we can use standard quasilinear elliptic
estimates~\cite{GilT83}, showing that a sequence $u_{R_{k}}$ converges locally in $C^{2}$ on $\wbar\Omega_{\epsilon}$ to
a $C^{2}$ function $u$ on $\wbar\Omega_{\epsilon}$ with $u=S_{\epsilon}$ on $\partial\Omega_{\epsilon}$, where
$S_{\epsilon}$ is defined relative to the original function $h$ and $u$ is a $C^{2,\alpha}$ solution of the SME.\end{proof}

\begin{proof}{\bf{}Completion of the proof of Theorem~\ref{ex-thm}:}
Let $S=\lim_{\epsilon\downarrow 0}S_{\epsilon}$; thus
\[%
S(r,y) = 
\psi(y)\tilde\varphi\bigl(r/\psi(y)\bigr)\text{ with } \psi(y) =e^{- 1/h(y)}, \quad r\ge 0,\,\, y\in U. 
\]%
If $u=u_{\epsilon}$ is the solution constructed in~\ref{periodic-ex} above, then, letting $\epsilon\downarrow 0$ and using
the estimates~$|Du_{\epsilon}|\le 2\alpha_{0}, \,|D_{y}u_{\epsilon}|\le \delta_{0}$ (true by construction of $u_{\epsilon}$),
we obtain (after taking a subsequence $\epsilon=\epsilon_{j}\downarrow 0$ if necessary) a non-negative Lipschitz function
$u$ (depending on $\tau$) on $\wbar\Omega$ with $|Du|\le 2\alpha_{0}$, $|D_{y}u|\le \delta_{0}$, $u=S$ on
$\partial\Omega$, $u$ a $C^{2,\alpha}$ solution of the SME on the open set where $u>0$, and
\[%
0\le u-\alpha_{0}|x|\le S-\alpha_{0}|x|\le C h^{j}(y) \text{ on } \Omega,\quad C=C(\ell,m,n,j).
\]%
It just remains to prove $u>\alpha_{0}|x|$ everywhere. But $u\ge \alpha_{0}|x|$ by construction, so if equality holds at some
point $(x_{0},y_{0})\in\Omega$ with $x_{0}\neq 0$ then the maximum principle (as in the discussion
following~\ref{lin-diff}) implies $u-\alpha|x|$ is identically zero in the connected component of $\Omega$ containing
$(x_{0},y_{0})$, contradicting the fact that $u=S>\alpha_{0}|x|$ on $\partial\Omega\setminus(\{0\}\times\R^{\ell})$.  On the
other hand if equality holds at $(0,y_{0})\in \Omega$ then we would have $u_{\epsilon_{j}}(0,y_{0})\to 0$, and by
Lemma~\ref{il-sw-app}\,(i) we would then have $u_{\epsilon_{j}}-\alpha_{0}|x|\to 0$ uniformly in some ball
$B_{\rho}(0,y_{0})$, so $u=\alpha_{0}|x|$ in $B_{\rho}(0,y_{0})$, and again we could conclude $u-\alpha_{0}|x|$ is
identically zero in a connected component of $\Omega$, contradicting the fact that $u>\alpha_{0}|x|$ on
$\partial\Omega\setminus(\{0\}\times\R^{\ell})$.\end{proof}

\newpage

\bibliographystyle{amsalpha}


\begin{thebibliography}{Law72}

\bibitem[FS20]{FouS20}
K.\ Fouladgar and L.~Simon, \emph{The symmetric minimal surface equation},
  Indiana Univ. Math. J. \textbf{69} (2020), 331--366.

\bibitem[GT83]{GilT83}
D.~Gilbarg and N.~Trudinger, \emph{Elliptic partial differential equations of
  second order}, 2nd ed., Springer-Verlag, Berlin-Heidelberg-New York, 1983.

\bibitem[HS85]{HarS85}
R.\ Hardt and L.~Simon, \emph{Area minimizing hypersurfaces with isolated
  singularities}, J.\ Reine u.\ Angew.\ Math. \textbf{362} (1985), 102--129.

\bibitem[Ilm96]{Ilm96}
T.~Ilmanen, \emph{A strong maximum principle for singular minimal
  hypersurfaces}, Calc. Var. and PDE \textbf{4} (1996), 443--467.

\bibitem[Law72]{Law72}
H.B. Lawson, \emph{The equivariant {P}lateau problem and interior regularity},
  Trans.\ Amer.\ Math.\ Soc. \textbf{173} (1972), 231--249.

\bibitem[Liu20]{Liu20}
Zhenhua Liu, \emph{Singularities of calibrated minimal surfaces can start and
  stop}, Preprint (2020).

\bibitem[Sim76]{Sim76}
L.~Simon, \emph{Interior gradient bounds for non-uniformly elliptic equations},
  Indiana Univ.\ Math.\ J. \textbf{25} (1976), 821--855.

\bibitem[Sim95]{Sim95b}
\bysame, \emph{Rectifiability of the singular sets of multiplicity 1 minimal
  surfaces and energy minimizing maps}, Surveys in Differential Geom.
  \textbf{II} (1995), 246--305.

\bibitem[Sim21]{Sim21a}
\bysame, \emph{A {L}iouville-type theorem for stable minimal hypersurfaces},
  Ars Inveniendi Analytica. arXiv:2101.06404v2 Paper No.5 (2021),
  35pp.

\bibitem[SW89]{SolW89}
B.~Solomon and B.~White, \emph{A strong maximum principle for varifolds that
  are stationary with respect to even parametric functionals}, Indiana Univ.\
  Math.\ J.\ \textbf{38} (1989), 683--691.

\end{thebibliography}

\newcommand{\noopsort}[1]{}
\providecommand{\bysame}{\leavevmode\hbox to3em{\hrulefill}\thinspace}
\providecommand{\MR}{\relax\ifhmode\unskip\space\fi MR }
\providecommand{\MRhref}[2]{%
  \href{http://www.ams.org/mathscinet-getitem?mr=#1}{#2}
}
\providecommand{\href}[2]{#2}

\enlargethispage*{1cm}

\medskip

\vbox{\parskip0pt\prevdepth-1000pt \multiply \baselineskip by 5 \divide \baselineskip by 6 \parskip0pt
  Mathematics Department  \\
  Stanford University\\
Stanford,
CA 94305}

lsimon@stanford.edu

\end{document}